\documentclass[11pt]{article}
\usepackage[T1]{fontenc}
\usepackage[utf8]{inputenc}
\usepackage{lmodern}
\usepackage[a4paper,margin=25mm]{geometry}
\usepackage{amsmath,amssymb,amsthm,mathtools,bm}
\usepackage{booktabs,array,enumitem}
\usepackage{microtype}
\usepackage[colorlinks=true,linkcolor=blue,urlcolor=blue,citecolor=blue]{hyperref}
\usepackage{fancyhdr}
\usepackage{float}
\usepackage{tikz}
\usetikzlibrary{decorations.markings}
\newsavebox\fitfigbox
\newlength\fitfigavail
\newlength\fitfignat
\newcommand\fitfigscale{1}
\makeatletter
\newcommand\fitfigrecord[1]{%
  \immediate\write\@auxout{\gdef\string\fitfigscale{#1}}}
\makeatother
\newcommand\fitfigure[2]{%
  \par
  \sbox\fitfigbox{#1}%
  \setlength\fitfignat{\dimexpr\ht\fitfigbox+\dp\fitfigbox\relax}%
  \ifdim\pagetotal=0pt
    \setlength\fitfigavail{\fitfignat}%
  \else
    \setlength\fitfigavail{\dimexpr\pagegoal-\pagetotal-106pt\relax}%
  \fi
  \pgfmathsetmacro\fitfigmax{min(1,\fitfigavail/\fitfignat)}%
  \pgfmathsetmacro\fitfiguse{min(\fitfigscale,\fitfigmax)}%
  \pgfmathsetmacro\fitfignew{floor(50*(\fitfigscale+\fitfigmax))/100}%
  \ifdim\fitfigmax pt<0.8pt
    \pgfmathsetmacro\fitfignew{floor(50*(\fitfigscale+1))/100}%
    \fitfigrecord{\fitfignew}%
    \begin{figure}[!htb]\centering
    \scalebox{\fitfigscale}{\usebox\fitfigbox}#2\end{figure}%
  \else
    \fitfigrecord{\fitfignew}%
    \begin{figure}[H]\centering
    \scalebox{\fitfiguse}{\usebox\fitfigbox}#2\end{figure}%
  \fi}
\hypersetup{pdftitle={Optimal long-time derivative estimates for underdamped Langevin dynamics}}
\allowdisplaybreaks[2]
\numberwithin{equation}{section}
\newtheorem{theorem}{Theorem}[section]
\newtheorem{proposition}[theorem]{Proposition}
\newtheorem{lemma}[theorem]{Lemma}
\newtheorem{corollary}[theorem]{Corollary}
\theoremstyle{definition}
\newtheorem{assumption}[theorem]{Assumption}
\newtheorem{remark}[theorem]{Remark}
\newcommand{\R}{\mathbb R}
\newcommand{\A}{\mathcal A}
\newcommand{\id}{I}
\newcommand{\dd}{\,\mathrm d}
\newcommand{\ip}[2]{\langle #1,#2\rangle_\mu}
\newcommand{\nm}[1]{\lVert #1\rVert_\mu}
\newcommand{\opnm}[1]{\lVert #1\rVert_{2\to2}}
\newcommand{\J}{\mathcal J}
\newcommand{\Sx}{\mathcal S_x}
\newcommand{\Lham}{\mathcal L_{\mathrm{ham}}}
\newcommand{\LFD}{\mathcal L_{\mathrm{FD}}}
\newcommand{\Cc}{C_c^\infty(\R^{2d})}
\newcommand{\cF}{\mathcal F}
\newcommand{\cG}{\mathcal G}
\newcommand{\cD}{\mathcal D}
\newcommand{\Kin}{\mathsf K}
\newcommand{\Vis}{\mathsf V}
\DeclareMathOperator{\Ree}{Re}
\DeclareMathOperator{\Imm}{Im}
\DeclareMathOperator{\Dom}{D}
\title{Optimal long-time derivative estimates for underdamped Langevin dynamics}
\author{LUKANG SUN\\[2pt]
{\normalsize Department of Mathematics, National University of Singapore}\\
{\normalsize\href{mailto:lukang.sun@nus.edu.sg}{\texttt{lukang.sun@nus.edu.sg}}}}
\date{November 14, 2026}

\begin{document}
\maketitle
\begin{abstract}
We study the long-time behaviour of derivatives of the semigroup
$P_t=e^{t\A}$ generated by the underdamped Langevin (kinetic
Fokker--Planck) operator $\A$ on $L^2(\mu)$, where the Gibbs measure $\mu$
is not required to be finite. Under a weighted relative bound on the
Hessian of the potential, we prove that for every integer $m\ge0$ the squared $L^2(\mu)$ norms of
$\A^mP_th_0$ and of $\nabla_{x,v}\A^mP_th_0$ are bounded by
$C_mt^{-2m}\|h_0\|^2$ and $D_mt^{-(2m+1)}\|h_0\|^2$ respectively for large
times, and the free dynamics shows that these exponents are optimal. The
proof rests on a resolvent estimate near the imaginary axis, valid for
every smooth potential, which follows from the adjoint relation between
the generator and the velocity reversal and from the Gaussian Poincar\'e
inequality in the velocity variable.
\end{abstract}

\noindent\textbf{Keywords:} underdamped Langevin dynamics, kinetic
Fokker--Planck equation, long-time behaviour, decay of derivatives,
resolvent estimates, operator semigroups, hypocoercivity, infinite Gibbs
measure

\medskip
\noindent\textbf{2020 Mathematics Subject Classification.} Primary 35Q84,
47D06; Secondary 35B40, 35H10, 47A10, 60J60, 82C31.

\tableofcontents

\section{Introduction}\label{sec:intro}

\subsection{The equation and its functional setting}\label{sec:setting}

Let $d\ge1$ be an integer, let $\gamma,\sigma>0$ be two fixed parameters,
and let $L\in C^\infty(\R^d)$ be a real-valued potential. We study the
kinetic Fokker--Planck equation
\begin{equation}\label{eq:A}
 \partial_th_t=\A h_t,\qquad
 \A=-\Lham+\gamma\LFD,
\end{equation}
for an unknown function $h_t=h_t(x,v)$ on the phase space
$\R^{2d}=\R^d_x\times\R^d_v$, where
\begin{equation}\label{eq:ops}
 \Lham=v\cdot\nabla_x-\nabla L\cdot\nabla_v,
 \qquad
 \LFD=\sigma\Delta_v-v\cdot\nabla_v .
\end{equation}
The operator $\Lham$ is the Liouville operator of the Hamiltonian flow with
energy $L(x)+\tfrac12|v|^2$, whereas $\LFD$ is the Ornstein--Uhlenbeck
operator in the velocity variable, with temperature $\sigma$. The
subscript stands for fluctuation--dissipation: the term $\sigma\Delta_v$
comes from the random force, which is the fluctuation, and the term
$-v\cdot\nabla_v$ from the friction, which is the dissipation. The
friction coefficient $\gamma$ is the factor with which $\LFD$ enters
$\A$. Equation \eqref{eq:A} describes the
underdamped Langevin dynamics
\begin{equation}\label{eq:SDE}
 \dd X_t=V_t\dd t,\qquad
 \dd V_t=-\nabla L(X_t)\dd t-\gamma V_t\dd t+\sqrt{2\gamma\sigma}\dd W_t
\end{equation}
in the following sense: if the law of $(X_t,V_t)$ has a density $f_t$ with
respect to the Lebesgue measure, then the density of this law relative to the
Gibbs measure
\begin{equation}\label{eq:mu}
 \dd\mu=e^{-(L(x)+|v|^2/2)/\sigma}\dd x\dd v ,
\end{equation}
namely $h_t=f_te^{(L+|v|^2/2)/\sigma}$, solves \eqref{eq:A}.
We emphasize that $\mu$ is \emph{not} required to be a finite measure in
this paper, and that no confinement is imposed on the position variable,
that is, we do not assume that $L(x)\to+\infty$ as $|x|\to\infty$. The only condition which we impose on $L$
concerns its Hessian, it is stated in Assumption~\ref{ass:WH}.

An integration by parts with respect to $\mu$ shows that $\Lham$ is
antisymmetric and that $\LFD$ is symmetric and nonpositive on test
functions; see \eqref{eq:ibp}. Consequently a solution of \eqref{eq:A}
satisfies, at least formally,
\[
 \frac{\dd}{\dd t}\int_{\R^{2d}}|h_t|^2\dd\mu
 =-2\gamma\sigma\int_{\R^{2d}}|\nabla_vh_t|^2\dd\mu .
\]
Thus the norm of the solution in $L^2(\mu)$ decreases, but the rate at
which it decreases involves only the derivatives of $h_t$ in the velocity
variable: the equation contains no term which directly damps the
variations of $h_t$ in the position variable.

The question addressed in this paper is how fast the derivatives of the
solution $h_t$ decay as $t\to\infty$ when nothing is assumed on the initial
datum beyond $h_0\in L^2(\mu)$, and when no functional inequality is
available in the position variable. The derivatives we consider are the
time derivatives of the solution, which are the powers of the generator
applied to it, $\partial_t^mh_t=\A^mh_t$, and the gradients of these
functions in the phase-space variables $x$ and $v$.

This question is motivated by the recent analysis of the overdamped
Langevin dynamics
\[
 \dd w_t=-\nabla L(w_t)\dd t+\sqrt{2\sigma}\dd B_t
\]
by Fornasier, Sun and Ward \cite{FSW}. Let $\pi=e^{-L/\sigma}\dd w$ denote
the Gibbs measure of that dynamics, which again is not assumed to be
integrable, let $\phi_t$ be the density of the law of $w_t$ relative to
$\pi$, and let $\mathfrak L=\sigma\Delta-\nabla L\cdot\nabla$ be the
generator. It is shown in \cite{FSW} that the Dirichlet energy
$\int|\nabla\phi_t|^2\dd\pi$ decays at the rate $O(1/t)$ for initial data
in $L^2(\pi)$, and more generally, for sufficiently regular potentials, that
the squared $L^2(\pi)$ norms of $\mathfrak L^m\phi_t$ and of
$\nabla\mathfrak L^m\phi_t$ decay like $t^{-2m}$ and $t^{-(2m+1)}$
respectively. The proof given in \cite{FSW} does not rely on spectral
theory, it uses only the symmetry of $\mathfrak L$ in $L^2(\pi)$: an
integration by parts shows that, along the flow, the time derivative of
each of the quantities
\[
 \int|\phi_t|^2\dd\pi,\qquad\int|\nabla\phi_t|^2\dd\pi,\qquad
 \int|\mathfrak L\phi_t|^2\dd\pi,\qquad
 \int|\nabla\mathfrak L\phi_t|^2\dd\pi,\qquad\dots
\]
is a negative multiple of the next one. All these quantities are therefore
nonincreasing in time, and since a nonincreasing function is bounded by
its average over a time interval, the decay of each of them follows from
that of the preceding one. The generator $\mathfrak L$ of the
overdamped dynamics is elliptic, and it is symmetric in $L^2(\pi)$. The
kinetic operator $\A$ has neither property: it is not symmetric in
$L^2(\mu)$, and its dissipation acts on the velocity variable only. It is
therefore not clear a priori that the same hierarchy of large-time
estimates survives for the underdamped dynamics, and the purpose of this
paper is to show that it does, with exactly the same powers of $t$.

\subsection{Assumption and main results}\label{sec:assumptions}

We work throughout in the complex Hilbert space
$L^2(\mu)=L^2(\R^{2d},\mu;\mathbb C)$, whose inner product and norm are
denoted by
\[
 \ip fg=\int_{\R^{2d}}f\,\overline g\dd\mu,\qquad \nm f^2=\ip ff .
\]
Complex-valued functions are needed
because the proof uses the resolvent $(z-\A)^{-1}$ at complex points $z$
close to the imaginary axis. This is no restriction: the operator $\A$ has
real coefficients, so that the semigroup which it generates maps
real-valued functions to real-valued functions, and all our results apply
in particular to real initial data. The only quantitative hypothesis of
the paper bounds the Hessian of $L$ relative to the weighted Dirichlet
form in the position variable.

\begin{assumption}[Weighted relative Hessian bound]\label{ass:WH}
There is a constant $C_L\ge1$ such that every $g\in\Cc$ satisfies
\begin{equation}\label{eq:WH}
 \int_{\R^{2d}}\|\nabla^2L(x)\|_{\mathrm{op}}^2\,|g|^2\dd\mu
 \le C_L\big(\nm g^2+\nm{\nabla_xg}^2\big).
\end{equation}
\end{assumption}

This assumption is taken from Villani's study of the kinetic
Fokker--Planck equation: it is the relative boundedness
condition~(7.2) of \cite{Villani}, and it also appears, in an equivalent
form, as \cite[Assumption~2]{Zhang}. It holds as soon as the Hessian is
controlled pointwise by the gradient in the sense that
\begin{equation}\label{eq:VG}
 \|\nabla^2L(x)\|_{\rm op}\le C\big(1+|\nabla L(x)|\big),\qquad x\in\R^d,
\end{equation}
for some constant $C\ge0$, by \cite[Lemma~A.18]{Villani}; this covers
potentials with bounded Hessian as well as potentials of arbitrary
polynomial growth. We stress that
Assumption~\ref{ass:WH} is an upper bound on the second derivatives of $L$
only: it implies neither confinement, nor a Poincar\'e inequality for the
measure $e^{-L/\sigma}\dd x$, nor the finiteness of this measure, and $L$
need not be bounded from below.

We denote by $P_t=e^{t\A}$ the semigroup of contractions on
$L^2(\mu)$ generated by $\A$, so that $h_t=P_th_0$, and by $\Dom(\A^m)$
the domain of the $m$th power of the generator; $\nabla_{x,v}=(\nabla_x,\nabla_v)$ is the
full gradient on $\R^{2d}$, and $H^1(\mu)$ is the weighted Sobolev space of
the functions $h\in L^2(\mu)$ whose weak gradient $\nabla_{x,v}h$ belongs
to $L^2(\mu)$. The precise details are given in
Section~\ref{sec:generator}. With these preparations, the main result
can be stated as follows.

\begin{theorem}[Long-time estimates for generator powers and for the gradient]
\label{thm:main}
Let $L\in C^\infty(\R^d)$ satisfy
Assumption~\ref{ass:WH}. Then, for every $t>0$ and every integer $m\ge0$,
\[
 P_tL^2(\mu)\subset\Dom(\A^m)
 \qquad\text{and}\qquad
 \A^mP_tL^2(\mu)\subset H^1(\mu).
\]
Moreover, for every integer $m\ge0$ there are finite constants $C_m$ and
$D_m$ such that every $h_0\in L^2(\mu)$ satisfies
\begin{align}
 \nm{\A^mP_th_0}^2&\le C_mt^{-2m}\nm{h_0}^2,
       &&t\ge1,\label{eq:main-even}\\
 \nm{\nabla_{x,v}\A^mP_th_0}^2
       &\le D_mt^{-(2m+1)}\nm{h_0}^2,
       &&t\ge2.\label{eq:main-odd}
\end{align}
Here $C_0=1$, and the other constants depend only on
$m,\gamma,\sigma$ and on the constant $C_L$ of
Assumption~\ref{ass:WH}, their values are given in
Table~\ref{tab:constants}.
\end{theorem}

Taking $m=0$ in \eqref{eq:main-odd} we see in particular that the full
gradient of the solution, including its position gradient, decays like
$t^{-1/2}$ in $L^2(\mu)$, uniformly over the unit ball in $L^2(\mu)$ of initial data.
More generally, the estimates \eqref{eq:main-even} and \eqref{eq:main-odd}
are the kinetic counterparts of the overdamped estimates of \cite{FSW}
recalled in Section~\ref{sec:setting}, with the same powers of $t$. The
thresholds $t\ge1$ and $t\ge2$
have no particular significance: for small times the regularization
estimates of Section~\ref{sec:energy} apply instead and give the usual
hypoelliptic short-time rates.

The proof of Theorem~\ref{thm:main} is given in Section~\ref{sec:spectral}.
It consists of two parts of a different nature. The first part is a bound on the
resolvent of the generator, and it holds for every smooth potential $L$,
without any further assumption. The generator
is related to its adjoint by the reversal of the velocity, and this
relation, combined with the Gaussian Poincar\'e inequality in the velocity
variable, gives this bound near the imaginary axis. It implies the decay of the powers of the
generator stated in the theorem, but only for initial data which are
sufficiently regular. The second part removes this restriction, and it is
only here that Assumption~\ref{ass:WH} enters: the energy method of
H\'erau \cite{Herau}, Villani \cite{Villani} and Zhang \cite{Zhang},
applied on a time interval of length one, shows that the semigroup
regularizes the initial datum, and it also controls the gradient of the
solution by the decrease of its norm, which gives the gradient estimate.

\begin{remark}[The semigroup and the law of the underdamped Langevin dynamics]
\label{rem:law}
Theorem~\ref{thm:main} is a statement about the semigroup $P_t$, that is,
about the solution $h_t=P_th_0$ of \eqref{eq:A} in $L^2(\mu)$. It is not, by
itself, a statement about the law of the underdamped Langevin dynamics
\eqref{eq:SDE}. The two are related by a result of Conrad and
Grothaus \cite[Corollary~2.3, Theorem~2.5, Remark~2.6 and
Corollary~2.7]{ConradGrothaus},
which holds for potentials that are bounded from below and locally
Lipschitz continuous; the second condition is satisfied by $C^\infty$
potential in this work. Under these conditions the underdamped Langevin dynamics
\eqref{eq:SDE} has a weak solution with infinite lifetime, and its
transition semigroup is the semigroup generated in $L^2(\mu)$ by the
closure of its generator, by Proposition~\ref{prop:adjoint} this semigroup
is the adjoint semigroup $P_t^*$. Consequently, if $L$ is bounded from below in addition to satisfying
Assumption~\ref{ass:WH}, and if the law of $(X_0,V_0)$ has a density
$h_0\in L^2(\mu)$ with respect to $\mu$, then for every $t>0$ the law of $(X_t,V_t)$ has the
density $P_th_0$ with respect to $\mu$, and Theorem~\ref{thm:main}
describes the decay of the derivatives of this relative density. When $L$
is not bounded from below, the underdamped Langevin dynamics
\eqref{eq:SDE} may explode in finite time, as it does for
$L(x)=-|x|^4$, whereas the semigroup $P_t$ is defined for every smooth
potential; in that case Theorem~\ref{thm:main} remains a statement about
$P_t$ only.
\end{remark}

\begin{remark}[Dependence of the constants on the potential]
\label{rem:CL}
The constants $C_m$ and $D_m$ of Theorem~\ref{thm:main} cannot be chosen
independently of the potential. This is in contrast with the overdamped
estimates of \cite{FSW}, whose constants depend only on $m$ and $\sigma$.
This is shown by the following example, which was obtained with the
assistance of the AI system Claude. Let $\omega>\max(1,\gamma/2)$ and $L(x)=\omega^2|x|^2/2$. This
potential satisfies Assumption~\ref{ass:WH} with $C_L=\omega^4$. Let
\[
 c=-\frac\gamma2+i\sqrt{\omega^2-\frac{\gamma^2}4}\,,
 \qquad\text{so that}\qquad
 c+\overline c=-\gamma,\qquad c\,\overline c=\omega^2 .
\]
Let $h(x,v)=v_1+c\,x_1$. Since
$\A v_1=\omega^2x_1-\gamma v_1$ and $\A x_1=-v_1$, we have
\[
 \A h=-(\gamma+c)\,v_1+\omega^2x_1
     =\overline c\,v_1+\overline c\,c\,x_1=\overline c\,h .
\]
Hence $h\in\Dom(\A)$ by \eqref{eq:maxdomain}, and
$P_th=e^{\overline ct}h$. Under the Gaussian measure $\mu$ the variables
$x_1$ and $v_1$ are independent, with variances $\sigma/\omega^2$ and
$\sigma$, so that $\nm h^2=2\sigma\,\mu(\R^{2d})$, whereas
$|\nabla_{x,v}h|^2=1+|c|^2=1+\omega^2$ is constant. Thus
$\nm{\nabla_{x,v}h}^2=(1+\omega^2)\nm h^2/(2\sigma)$, and therefore
\[
 \nm{\A P_th}^2=\omega^2e^{-\gamma t}\nm h^2,\qquad
 \nm{\nabla_{x,v}P_th}^2=\frac{1+\omega^2}{2\sigma}\,e^{-\gamma t}\nm h^2 .
\]
Comparing with \eqref{eq:main-even} for $m=1$ at $t=1$, and with
\eqref{eq:main-odd} for $m=0$ at $t=2$, we obtain
\[
 C_1\ge e^{-\gamma}\omega^2=e^{-\gamma}C_L^{1/2},\qquad
 D_0\ge\frac{e^{-2\gamma}}{\sigma}\,(1+\omega^2)
      =\frac{e^{-2\gamma}}{\sigma}\,\big(1+C_L^{1/2}\big) .
\]
The reason is that $\A$ is not selfadjoint, its eigenvalue $\overline c$
has the real part $-\gamma/2$, whatever the value of $\omega$, and an
imaginary part of size $\omega$, which corresponds to fast oscillations of
the solution.
\end{remark}

\begin{remark}[Potentials of class $C^2$]\label{rem:C2}
Theorem~\ref{thm:main} extends to potentials which are only twice
continuously differentiable. More precisely, let $L\in C^2(\R^d)$ be
bounded from below and satisfy the pointwise condition \eqref{eq:VG}. Then
the closure of $\A$ on $\Cc$ generates a contraction semigroup $P_t$ on
$L^2(\mu)$ \cite[Corollary~2.3 and Remark~2.4]{ConradGrothaus}, and the
conclusions of
Theorem~\ref{thm:main} hold for this semigroup, with constants which
depend only on $m,d,\gamma,\sigma$ and on the constant $C$ of
\eqref{eq:VG}. The proof is given in Appendix~\ref{apx:C2}.

Here $\A$ denotes this closure, and its powers are defined by iteration, as
in Section~\ref{sec:generator}: $\Dom(\A^m)$ is the set of the functions
$h\in\Dom(\A)$ such that $\A h,\dots,\A^{m-1}h$ also belong to $\Dom(\A)$,
and $\A^mh=\A(\A^{m-1}h)$. When $L$ is smooth, $\A^mh$ can also be computed
by expanding the differential expression $\A^m$ and applying it to $h$ in
the sense of distributions. When $L$ is only twice continuously
differentiable, this is no longer possible: for $m\ge3$ the expanded
expression involves derivatives of $L$ of order higher than two, and it has
no classical meaning. This does not mean that $\Dom(\A^m)$ is trivial. For
example, let $d=1$ and $L(x)=|x|^{5/2}$. This potential is bounded from
below, satisfies \eqref{eq:VG} and is twice continuously differentiable,
but its third derivative is unbounded near the origin. Since $L$ is smooth
away from the origin, a test function $\psi(x,v)$ which vanishes for $x$
in a neighbourhood of the origin is mapped by $\A$ to a test function of
the same kind, and therefore $\psi\in\Dom(\A^m)$ for every $m$. In the same
way, since $P_th_0\in\Dom(\A^m)$ by the first conclusion of
Theorem~\ref{thm:main}, the function $\A^mP_th_0$ which appears in
\eqref{eq:main-even} and \eqref{eq:main-odd} is well defined by the
iterative definition.
\end{remark}

The powers of $t$ in Theorem~\ref{thm:main} cannot be improved under
Assumption~\ref{ass:WH}, as the following explicit analysis of the free
dynamics shows; in its statement
$\opnm T=\sup_{h\ne0}\nm{Th}/\nm h$ denotes the operator norm of a bounded
operator $T$ on $L^2(\mu)$.

\begin{theorem}[Optimality of the exponents]\label{thm:sharp}
Let $L=0$, which satisfies Assumption~\ref{ass:WH} with $C_L=1$. Then, for
all $t>0$,
\begin{align}
 \opnm{\A^mP_t}^2&\ge\left(\frac m{et}\right)^{2m},
     &&m\ge1,\label{eq:sharp-even}\\
 \sup_{h\ne0}\frac{\nm{\nabla_{x,v}\A^mP_th}^2}{\nm h^2}
   &\ge\frac\gamma\sigma\big(1+\gamma^{-2}\big)
        \left(\frac{2m+1}{2et}\right)^{2m+1},
     &&m\ge0.\label{eq:sharp-odd}
\end{align}
Consequently the powers of $t$ in \eqref{eq:main-even} and
\eqref{eq:main-odd} cannot be improved in any class of potentials containing
$L=0$.
\end{theorem}

The proof of Theorem~\ref{thm:sharp} is given in Section~\ref{sec:sharp}.
The lower bounds are obtained by exhibiting a closed subspace of
$L^2(\mu)$ which is invariant under the free dynamics.

\subsection{Related literature}\label{sec:literature}

This paper is related to two bodies of work: the abstract theory of the
decay rates of operator semigroups, to which our method belongs, and the
literature on the long-time behaviour of the kinetic Fokker--Planck
equation, with which Theorem~\ref{thm:main} is to be compared.

\paragraph{Decay rates of operator semigroups.}
There is a general theory which deduces the rate of decay of a bounded
semigroup on a Hilbert space from the size of the resolvent of its
generator at the points $is$ of the imaginary axis. If the resolvent is
bounded on the whole imaginary axis, the semigroup decays exponentially;
this is the theorem of Gearhart, Pr\"uss and Greiner
\cite[Theorem~V.1.11]{EngelNagel}. If the resolvent exists on the whole
imaginary axis but grows like a power of $|s|$ as $|s|\to\infty$, then
the solutions whose initial datum belongs to the domain of the generator
decay polynomially, at a rate which is determined by this power. This is
due to Batty and Duyckaerts \cite{BattyDuyckaerts} and, with the optimal
rates, to Borichev and Tomilov \cite{BorichevTomilov}; further
developments are contained in \cite{BCT} and \cite{RSS}, and the theory is
surveyed in \cite{CST}. Finally, if the origin belongs to the spectrum of
the generator, the resolvent is unbounded as $s\to0$. For this case Batty,
Chill and Tomilov \cite{BCT} showed that a bound of the resolvent by a
multiple of $1/|s|$ near the origin, together with a uniform bound for
large $|s|$, is equivalent to the following property: the time derivative
of the solution decays like $1/t$ whenever the initial datum belongs to
the domain of the generator.

The operator $\A$ of the present paper is in the last situation, because
the origin belongs in general to its spectrum when the measure $\mu$ is
infinite. From the point of view of this theory, what we do is the
following. We prove the required bound on the resolvent of $\A$, for every
smooth potential (Proposition~\ref{prop:res}). The theory then gives the
decay of the time derivatives of the solution for initial data in the
domain of a power of $\A$ (Proposition~\ref{prop:reg}). Finally we show
that the semigroup maps $L^2(\mu)$ into these domains at every positive
time (Proposition~\ref{prop:smooth}), so that the decay holds for all
initial data in $L^2(\mu)$.

\paragraph{Long-time behaviour of the underdamped dynamics.}
The existing long-time results for the kinetic Fokker--Planck equation
either concern the solution itself or, when they also control its
derivatives, rely on a spectral gap. Under a functional inequality in the
position variable and for a finite Gibbs measure, H\'erau \cite{Herau} and
Villani \cite[Section~7]{Villani} prove exponential convergence to
equilibrium, also in norms which control the derivatives of the solution;
the relative bounds on the Hessian used by Villani are the model for
Assumption~\ref{ass:WH}, and explicit rates are given by Cao, Lu and Wang
\cite{CaoLuWang}. Without confinement, or when the Gibbs measure is not
integrable, the solution decays algebraically at a rate which depends on
the dimension, for initial data in spaces smaller than $L^2(\mu)$: this is
shown in \cite{BDMMS,BDS} by means of Nash-type inequalities, and in
\cite{Wang,PWZ}, together with the leading profile of the solution, by a
spectral analysis of the resolvent near the origin. For weakly confining
potentials, for which the Gibbs measure is finite but no spatial Poincar\'e
inequality holds, stretched-exponential or algebraic rates of convergence to
equilibrium are obtained in \cite{Cao,BSWW} under additional assumptions on
the initial datum; the weighted Poincar\'e inequalities for heavy-tailed
measures which are relevant there go back to \cite{BobkovLedoux}. In
contrast with all these works, Theorem~\ref{thm:main} concerns the
derivatives of the solution, requires neither a functional inequality nor a
finite Gibbs measure, and uses only the $L^2(\mu)$ norm of the initial
datum. The powers of $t$ which it gives and the constants $C_m$ and
$D_m$ do not depend on the dimension $d$. In return, the theorem gives
no information on the decay of the solution itself.

\subsection{Organization of the paper}

Section~\ref{sec:generator} collects the preliminaries on the generator
$\A$. Section~\ref{sec:spectral} contains the proof of
Theorem~\ref{thm:main}: the resolvent estimate and the regularized decay
estimate in Section~\ref{sec:resolvent}, the smoothing property in
Section~\ref{sec:energy}, and the conclusion in
Section~\ref{sec:proofmain}. Section~\ref{sec:sharp} contains the proof of
Theorem~\ref{thm:sharp}. The appendices contain the proofs of
Propositions~\ref{prop:maxdissipative} and~\ref{prop:adjoint}
(Appendix~\ref{apx:maxdissipative}), of Lemmas~\ref{lem:reverse}
and~\ref{lem:static} (Appendix~\ref{apx:energy}), of the elliptic
regularization on which the latter rely (Appendix~\ref{apx:elliptic}),
and of the extension to potentials of class $C^2$ stated in
Remark~\ref{rem:C2} (Appendix~\ref{apx:C2}).

\paragraph{Notation.}
Throughout, $C$ and $c$ denote positive constants whose value may change
from one line to the next. For an integer $k\ge1$, $H^k(\mu)$ denotes the
weighted Sobolev space of the functions in $L^2(\mu)$ whose weak
derivatives up to order $k$ belong to $L^2(\mu)$, and $H^1_x(\mu)$ and
$H^1_v(\mu)$ denote the spaces of the functions in $L^2(\mu)$ whose weak
gradient in $x$, respectively in $v$, belongs to $L^2(\mu)$. The test
functions $\Cc$ are dense in all these spaces.

\section{The generator, its adjoint, and velocity reversal}\label{sec:generator}

This section collects the functional-analytic facts on which the rest of
the paper relies. None of them requires an assumption on the Hessian of
$L$: the potential is only assumed to be smooth. These facts are known,
but they are scattered in the literature, where they are stated under
various assumptions and in different forms. We gather them here in the
form in which they will be used, and for the sake of completeness we give
the proofs of the two propositions in Appendix~\ref{apx:maxdissipative}.
The references are collected at the end of the section.

\paragraph{Integration by parts.}
We first record the two integrations by parts which express the structure
of $\A$ in $L^2(\mu)$. Since
$\LFD=\sigma e^{|v|^2/(2\sigma)}\nabla_v\cdot
(e^{-|v|^2/(2\sigma)}\nabla_v)$, and since the vector field
$(v,-\nabla L)$ defining $\Lham$ is divergence free and annihilates the
weight $e^{-(L+|v|^2/2)/\sigma}$, we have for all $f,g\in\Cc$
\begin{equation}\label{eq:ibp}
 \ip{\LFD f}{g}=-\sigma\ip{\nabla_vf}{\nabla_vg},\qquad
 \ip{\Lham f}{g}=-\ip f{\Lham g}.
\end{equation}
Consequently the formal adjoint of $\A=-\Lham+\gamma\LFD$ in $L^2(\mu)$ is
the operator $\Lham+\gamma\LFD$, in the sense that
\begin{equation}\label{eq:formaladjoint}
 \ip{\A f}{g}=\ip f{(\Lham+\gamma\LFD)g},\qquad f,g\in\Cc,
\end{equation}
and taking $g=f$ and real parts gives the dissipation identity on test
functions,
\begin{equation}\label{eq:dirichlet}
 -\Ree\ip{\A h}{h}=\gamma\sigma\nm{\nabla_vh}^2 ,\qquad h\in\Cc .
\end{equation}
The same integrations by parts give the formal adjoints in $L^2(\mu)$ of
the gradients $\nabla_v$ and $\nabla_x$. They map a vector field
$F=(F_1,\dots,F_d)$ with components in $\Cc$ to the functions
\begin{equation}\label{eq:adjoints}
 \nabla_v^*F=-\nabla_v\cdot F+\frac{v\cdot F}\sigma,\qquad
 \nabla_x^*F=-\nabla_x\cdot F+\frac{\nabla L\cdot F}\sigma ,
\end{equation}
that is, $\ip{\nabla_vf}{F}=\ip f{\nabla_v^*F}$ and
$\ip{\nabla_xf}{F}=\ip f{\nabla_x^*F}$ for $f\in\Cc$. We apply these
formulas to the vector fields $\nabla_xh$ and $\nabla_vh$, where
$h\in\Cc$. In the difference $\nabla_v^*\nabla_xh-\nabla_x^*\nabla_vh$ the
terms which contain the mixed derivatives cancel, and comparing with the
definition \eqref{eq:ops} of $\Lham$ and $\LFD$ we obtain
\begin{equation}\label{eq:adjointform}
 \Lham h=\sigma\big(\nabla_v^*\nabla_xh-\nabla_x^*\nabla_vh\big),\qquad
 \LFD h=-\sigma\,\nabla_v^*\nabla_vh .
\end{equation}

\paragraph{Notions from the theory of semigroups.}
We recall the few notions from the theory of operator semigroups which are
needed, following Engel and Nagel \cite[Section~II.3]{EngelNagel}. Let $T$
be a linear operator on $L^2(\mu)$ defined on a dense subspace $\Dom(T)$,
its domain. It is called \emph{dissipative} if
\[
 \Ree\ip{Th}{h}\le0,\qquad h\in\Dom(T).
\]
A densely defined dissipative operator is \emph{closable}
\cite[Proposition~II.3.14(iv)]{EngelNagel}, which means that
\begin{equation}\label{eq:closable}
 h_n\in\Dom(T),\quad h_n\longrightarrow0\quad\text{and}\quad
 Th_n\longrightarrow g\ \text{ in }L^2(\mu)
 \qquad\Longrightarrow\qquad g=0 ,
\end{equation}
and its closure is again dissipative. A dissipative operator $T$ is called
\emph{maximal dissipative} if, for some $\lambda>0$, the operator
$\lambda-T$ maps $\Dom(T)$ onto $L^2(\mu)$. By the Lumer--Phillips theorem
\cite[Theorem~II.3.15]{EngelNagel}, such an operator generates a strongly
continuous semigroup of contractions. Recall that the \emph{resolvent set}
$\rho(T)$ of a operator $T$ is the set of complex numbers $z$ for
which $z-T$ maps $\Dom(T)$ bijectively onto $L^2(\mu)$; the inverse
$(z-T)^{-1}$, called the resolvent of $T$ at $z$, is then a bounded
operator on $L^2(\mu)$ \cite[Section~IV.1]{EngelNagel}. For a maximal
dissipative operator $T$, every complex number $z$ with $\Ree z>0$ belongs
to $\rho(T)$, with
\begin{equation}\label{eq:halfplane}
 \opnm{(z-T)^{-1}}\le\frac1{\Ree z},\qquad\Ree z>0,
\end{equation}
by \cite[Theorem~II.3.5]{EngelNagel}. On a Hilbert space, a densely
defined dissipative operator is maximal dissipative if and only if it can not be extended to a strictly larger domain
\cite[Corollary to Theorem~1.1.1]{Phillips}.

The operator $\A$ defined on $\Cc$ is densely defined, because $\Cc$ is
dense in $L^2(\mu)$, and it is dissipative by \eqref{eq:dirichlet}. It is
therefore closable by \cite[Proposition~II.3.14(iv)]{EngelNagel}, which
states that a densely defined dissipative operator on a Banach space is
closable and that its closure is again dissipative. From now on $\A$
denotes its closure. Explicitly,
$\Dom(\A)$ consists of those $h\in L^2(\mu)$ for which there exist
$h_n\in\Cc$ and $g\in L^2(\mu)$ such that
\[
 h_n\longrightarrow h,\qquad \A h_n\longrightarrow g
                       \quad\text{in }L^2(\mu),
\]
and for such $h$ one sets $\A h=g$; applying \eqref{eq:closable} to the
difference of two approximating sequences shows that $g$ does not depend on
the sequence.

\paragraph{Maximal dissipativity and the operator domain.}
The first proposition describes the closure of $\A$ and the semigroup which
it generates.

\begin{proposition}[Maximal dissipativity and the operator domain]
\label{prop:maxdissipative}
Let $L\in C^\infty(\R^d)$ and let $\gamma,\sigma>0$. The
closure in $L^2(\mu)$ of $\A$, initially defined on $\Cc$, is maximal
dissipative, and its domain is the \emph{maximal domain}
\begin{equation}\label{eq:maxdomain}
 \Dom(\A)=\{h\in L^2(\mu):\A h\in L^2(\mu)
                                  \text{ in the sense of distributions}\}.
\end{equation}
It generates a strongly continuous semigroup of contractions
$P_t=e^{t\A}$ on $L^2(\mu)$. Thus $P_0=\id$, $P_{t+s}=P_tP_s$,
\[
 \nm{P_th}\le\nm h\qquad\text{and}\qquad
 \lim_{t\downarrow0}\nm{P_th-h}=0,\qquad h\in L^2(\mu),\quad t\ge0,
\]
and for every $h_0\in\Dom(\A)$ the function $h_t=P_th_0$ takes its values in
$\Dom(\A)$, is continuously differentiable with values in $L^2(\mu)$, and is
the unique solution of \eqref{eq:A} with initial datum $h_0$.
\end{proposition}

For an integer $m\ge1$ we denote by $\Dom(\A^m)$ the domain of the $m$th
power of $\A$, that is, the set of $h\in\Dom(\A)$ such that
$\A h,\dots,\A^{m-1}h$ also belong to $\Dom(\A)$.

\paragraph{The adjoint and velocity reversal.}
Define the velocity-reversal operator by
\[
 (\J h)(x,v)=h(x,-v).
\]
Since $\mu$ is invariant under $v\mapsto-v$, the operator $\J$ is a
selfadjoint unitary involution on $L^2(\mu)$, that is,
$\J^*=\J$, $\J^2=\id$ and $\nm{\J h}=\nm h$; the operator
$\tfrac12(\id-\J)$ is the orthogonal projection onto the functions which
are odd in the full velocity variable. On test functions $\J$ commutes with
$\LFD$ and anticommutes with $\Lham$, so that
\begin{equation}\label{eq:Jcore}
 \J\A\J=\Lham+\gamma\LFD\qquad\text{on }\Cc ,
\end{equation}
which by \eqref{eq:formaladjoint} is the formal adjoint of $\A$. So far
this is only an identity between differential operators acting on test
functions. The next proposition states that it remains true for the
closure of $\A$ and its adjoint $\A^*$: the operators $\A^*$ and $\J\A\J$
have the same domain, and they coincide on it.

\begin{proposition}[Adjoint and velocity reversal]\label{prop:adjoint}
Let $L\in C^\infty(\R^d)$. The adjoint of $\A$ in
$L^2(\mu)$ is the closure of $\Lham+\gamma\LFD$ on $\Cc$, and
\begin{equation}\label{eq:J}
 \A^*=\J\A\J,\qquad
 \Dom(\A^*)=\J\Dom(\A)
           =\{h:\J h\in\Dom(\A)\}.
\end{equation}
Consequently:
\begin{enumerate}[label=\textup{(\roman*)},leftmargin=*]
\item the operator $\J\A$ is selfadjoint on $\Dom(\A)$, and in particular
\begin{equation}\label{eq:Jreal}
 \Imm\ip{\A h}{\J h}=0,\qquad h\in\Dom(\A);
\end{equation}
\item every $h\in\Dom(\A)$ belongs to $H^1_v(\mu)$, and the dissipation
identity \eqref{eq:dirichlet} holds on all of $\Dom(\A)$.
\end{enumerate}
\end{proposition}

The selfadjointness of $\J\A$ can be expressed in terms of the pairing
\[
 \ip{\J h}{g}
       =\int h(x,-v)\,\overline{g(x,v)}\dd\mu .
\]
This pairing is Hermitian, which means that exchanging $h$ and $g$ replaces
it by its complex conjugate, $\ip{\J g}{h}=\overline{\ip{\J h}{g}}$; this
holds because $\J$ is selfadjoint, and it implies that $\ip{\J h}h$ is
real. The pairing is not positive, however, since $\ip{\J h}h=-\nm h^2$
when $h$ is odd in the velocity. The selfadjointness of $\J\A$ says that
$\A$ is symmetric for this pairing,
\[
 \ip{\J\A h}{g}=\ip{\J h}{\A g},
                   \qquad h,g\in\Dom(\A).
\]
The diagonal value of the pairing is
\begin{equation}\label{eq:Jform}
 \ip{\J h}h
 =\nm h^2-\tfrac12\nm{(\id-\J)h}^2 ,
\end{equation}
as one sees by expanding $\nm{(\id-\J)h}^2=2\nm h^2-2\ip{\J h}h$.

\paragraph{Dissipation on the odd velocity subspace.}
The second structural fact is a coercivity estimate for the velocity-odd
part of a function. It is supplied by the model itself, through the
Gaussian Poincar\'e inequality in the velocity variable, and it involves no
assumption on the position measure; in particular it remains valid when
$\mu$ is infinite.

\begin{lemma}[Dissipation on odd functions]\label{lem:odd}
For every $h\in\Dom(\A)$,
\begin{equation}\label{eq:odd}
 \tfrac14\nm{(\id-\J)h}^2
   \le\tfrac{\sigma}{4}\nm{\nabla_v(\id-\J)h}^2
   \le\sigma\nm{\nabla_vh}^2,
 \quad\text{and hence}\quad
 -\Ree\ip{\A h}h\ge\tfrac\gamma4\nm{(\id-\J)h}^2 .
\end{equation}
\end{lemma}
\begin{proof}
By Proposition~\ref{prop:adjoint}(ii) we have $h\in H^1_v(\mu)$, and the
same holds for $\J h$ because $\nabla_v(\J h)=-\J(\nabla_vh)$. For almost
every fixed $x$, the function $v\mapsto((\id-\J)h)(x,v)$ is odd and
therefore has zero mean with respect to the centred Gaussian measure with
density proportional to $e^{-|v|^2/(2\sigma)}$, whose Poincar\'e constant is
$\sigma$. The Gaussian Poincar\'e inequality thus gives, for almost every
$x$,
\[
 \int_{\R^d}|(\id-\J)h|^2e^{-|v|^2/(2\sigma)}\dd v
   \le\sigma\int_{\R^d}|\nabla_v(\id-\J)h|^2e^{-|v|^2/(2\sigma)}\dd v .
\]
Multiplying by $e^{-L(x)/\sigma}$ and integrating in $x$, which is
legitimate whether or not the Gibbs measure has finite mass, proves the
first inequality in \eqref{eq:odd}. For the second we write
$\nabla_v(\id-\J)h=(\id+\J)\nabla_vh$ and use that $\id+\J$ has norm at most
two. The last assertion follows from the dissipation identity
\eqref{eq:dirichlet} on $\Dom(\A)$.
\end{proof}

As we said at the beginning of this section, its statements are known in
one form or another, and we indicate here where they can be found. The maximal dissipativity asserted in
Proposition~\ref{prop:maxdissipative} is due to Helffer and Nier
\cite[Proposition~5.5]{HelfferNier}, who proved it for an arbitrary smooth
potential in the unweighted space $L^2(\dd x\dd v)$, to which $L^2(\mu)$ is
unitarily equivalent. In the weighted space $L^2(\mu)$ their result is
quoted by Villani \cite[Section~2.4]{Villani} for confining potentials, and
a statement formulated directly in weighted spaces, for locally Lipschitz
potentials, is \cite[Theorem~2.9]{GrothausNonnenmacher}; in both references
$\gamma=\sigma=1$. The description \eqref{eq:maxdomain} of the domain is
stated for the unweighted operator in \cite[Remark~3.2]{NierLectures}.

The adjoint relation \eqref{eq:J} is known as the generalized PT-symmetry
of kinetic Fokker--Planck operators. The letters P and T stand for parity
and time reversal: in quantum mechanics an operator is called PT-symmetric
if it is invariant under the reflection of the space variable combined
with the reversal of time. In the kinetic setting the reversal of time is
the reversal of the velocity, $v\mapsto-v$, and the word ``generalized''
indicates that this reversal does not leave the generator invariant but
exchanges it with its adjoint. For the unweighted operator this relation
is
noted in \cite[proof of Corollary~6.5]{HelfferNier} and in
\cite[Section~5.1]{NierLectures}, and it was exploited by H\'erau, Hitrik
and Sj\"ostrand \cite[Section~3.1, equations~(3.3)--(3.4)]{HHS}, who used
the pairing $\ip{\J h}{g}$ to show that the small
eigenvalues are real in a semiclassical multiple-well problem; see also
\cite[Section~3.2, Proposition~3.2]{Hitrik},
\cite[Section~15.3.2.2, equation~(15.5)]{NierChapter} and, for geometric
kinetic Fokker--Planck operators with boundary conditions,
\cite[Section~8.2]{NierBoundary}. In weighted spaces, the fact that the
adjoint of the generator is the closure of its formal adjoint, deduced from
the velocity reversal, is \cite[Corollary~2.12]{GrothausNonnenmacher}. For
confining potentials the same fact is proved by Breiten and Kunisch
\cite[Proposition~1]{BreitenKunisch}, together with the inclusion
$\Dom(\A)\subset H^1_v(\mu)$ and the validity of the dissipation identity on
all of $\Dom(\A)$ \cite[Lemma~3]{BreitenKunisch}. The equality of domains in
\eqref{eq:J} is stated explicitly, for
the unweighted operator with a potential whose Hessian is bounded, in
\cite[Proposition~A.10]{ABL}.

We could not find Propositions~\ref{prop:maxdissipative}
and~\ref{prop:adjoint} in the literature in exactly the present form, that
is, in $L^2(\mu)$, with general $\gamma$ and $\sigma$, and for an arbitrary
smooth potential. For this reason their proofs are given in
Appendix~\ref{apx:maxdissipative}.

\section{Proof of Theorem~\ref{thm:main}}
\label{sec:spectral}

This section contains the three ingredients of the proof of
Theorem~\ref{thm:main}. The first two, the resolvent bound and the
regularized decay estimate, use no assumption on the Hessian of $L$;
Assumption~\ref{ass:WH} enters only in the third one. They are combined
in Section~\ref{sec:proofmain}.

\subsection{The resolvent bound and the regularized decay estimate}
\label{sec:resolvent}\label{sec:contour}

We now combine the identity \eqref{eq:Jreal}, which is a consequence of the adjoint relation \eqref{eq:J}, with the coercivity
\eqref{eq:odd} to bound the resolvent of $\A$ in a strip to the left of the
imaginary axis, away from the real axis.

\begin{proposition}[Resolvent bound]\label{prop:res}
Let $L\in C^\infty(\R^d)$. Then
\[
 \{z\in\mathbb C:\Ree z\ge-\gamma/4,\ \Imm z\ne0\}\subset\rho(\A),
\]
and for every $z$ in this set
\begin{equation}\label{eq:Rbound}
 \opnm{(z-\A)^{-1}}\le\frac4\gamma+\frac2{|\Imm z|}.
\end{equation}
\end{proposition}

This resolvent bound is the key estimate of the paper. Its proof uses only
three properties of $\A$: its maximal dissipativity, the fact that
$\ip{\A h}{\J h}$ is real, which is \eqref{eq:Jreal}, and the coercivity
\eqref{eq:odd} on odd functions. In particular, no bound on the Hessian of
$L$ is needed.

\begin{proof}
Let $h\in\Dom(\A)$ and $z\in\mathbb C$. Taking the real part of the pairing
of $(z-\A)h$ with $h$ and using \eqref{eq:odd} gives
\begin{equation}\label{eq:real-res}
 \tfrac\gamma4\nm{(\id-\J)h}^2
 \le-\Ree z\,\nm h^2+\nm{(z-\A)h}\nm h.
\end{equation}
Pairing instead with $\J h$ and taking imaginary parts, the contribution of
$\A$ vanishes by \eqref{eq:Jreal}, while $\ip h{\J h}$ is real and is given
by \eqref{eq:Jform}; hence
\begin{equation}\label{eq:imag-res}
 \Imm z\left(\nm h^2-\tfrac12\nm{(\id-\J)h}^2\right)
       =\Imm\ip{(z-\A)h}{\J h}.
\end{equation}
Assume now that $\Imm z\ne0$. Dividing \eqref{eq:imag-res} by $|\Imm z|$,
bounding its right-hand side by $\nm{(z-\A)h}\nm h$, and inserting
\eqref{eq:real-res} to control $\tfrac12\nm{(\id-\J)h}^2$, we obtain
\begin{equation}\label{eq:apriori}
 \left(1+\frac{2\Ree z}{\gamma}\right)\nm h^2
 \le\left(\frac2\gamma+\frac1{|\Imm z|}\right)
                         \nm{(z-\A)h}\nm h.
\end{equation}
If $\Ree z\ge-\gamma/4$, the coefficient on the left-hand side is at least
$1/2$, and we arrive at the a priori bound
\[
 \nm h\le\left(\frac4\gamma+\frac2{|\Imm z|}\right)\nm{(z-\A)h},
 \qquad h\in\Dom(\A).
\]
This inequality shows that $z-\A$ is injective on $\Dom(\A)$: if
$(z-\A)h=0$, then $h=0$. It does not yet show that $z$ belongs to the
resolvent set, because for this $z-\A$ must also map $\Dom(\A)$ onto
$L^2(\mu)$. We prove the surjectivity by starting from the right half-plane,
where it is known, and moving to the left in small steps.

Fix a real number $\beta\ne0$, put
\[
 M=\frac4\gamma+\frac2{|\beta|},
\]
and consider the points $z_s=s+i\beta$ of the horizontal line
$\Imm z=\beta$, for $s\ge-\gamma/4$. All these points have the same
imaginary part, so that the a priori bound holds for all of them with the
same constant:
\begin{equation}\label{eq:apriori-line}
 \nm h\le M\,\nm{(z_s-\A)h},\qquad h\in\Dom(\A),\quad s\ge-\gamma/4 .
\end{equation}
The key observation is the following: if $z_{s_0}\in\rho(\A)$ for some
$s_0\ge-\gamma/4$, then $z_s\in\rho(\A)$ for every $s\ge-\gamma/4$ with
$|s-s_0|<1/M$. Indeed, let $f\in L^2(\mu)$ and apply
\eqref{eq:apriori-line} to $h=(z_{s_0}-\A)^{-1}f$; this gives
$\nm{(z_{s_0}-\A)^{-1}f}\le M\nm f$, that is,
$\opnm{(z_{s_0}-\A)^{-1}}\le M$. Since $z_s-z_{s_0}=s-s_0$, we can write
\[
 z_s-\A=(z_{s_0}-\A)+(s-s_0)
       =\Big[\id+(s-s_0)(z_{s_0}-\A)^{-1}\Big](z_{s_0}-\A)
 \qquad\text{on }\Dom(\A).
\]
The operator $B=(s-s_0)(z_{s_0}-\A)^{-1}$ is bounded on $L^2(\mu)$ with
$\opnm B\le|s-s_0|M<1$, and therefore $\id+B$ is invertible on $L^2(\mu)$,
its inverse being the sum of the convergent series
$\sum_{n\ge0}(-B)^n$. Thus $z_s-\A$ is the composition of the bijection
$z_{s_0}-\A$ from $\Dom(\A)$ onto $L^2(\mu)$ with the bijection $\id+B$ of
$L^2(\mu)$ onto itself. It is therefore a bijection from $\Dom(\A)$ onto
$L^2(\mu)$, which means that $z_s\in\rho(\A)$.

It remains to apply this observation repeatedly. For $s>0$ the point $z_s$
belongs to $\rho(\A)$, because $\A$ is maximal dissipative, as recalled in
\eqref{eq:halfplane}. Starting from $s_0=1/(2M)$, the observation shows that
$z_s\in\rho(\A)$ for all $s\ge-\gamma/4$ with $|s-s_0|<1/M$, in particular
for $s=0$; starting again from $s_0=0$ we reach $s=-1/(2M)$, and so on.
Each step moves to the left by $1/(2M)$, a length which does not depend on
the starting point because the constant in \eqref{eq:apriori-line} does not
depend on $s$. After finitely many steps we have covered the whole interval
$[-\gamma/4,0]$. Hence $z_s\in\rho(\A)$ for every $s\ge-\gamma/4$, and since
$\beta\ne0$ was arbitrary, every $z$ with $\Ree z\ge-\gamma/4$ and
$\Imm z\ne0$ belongs to $\rho(\A)$. Finally, for such $z$ and for
$f\in L^2(\mu)$, the a priori bound applied to $h=(z-\A)^{-1}f$ is exactly
\eqref{eq:Rbound}.
\end{proof}

The estimates \eqref{eq:real-res} and \eqref{eq:imag-res} have a simple
interpretation. Every function $h$ is the sum of its even part
$\tfrac12(\id+\J)h$ and of its odd part $\tfrac12(\id-\J)h$ in the velocity
variable, and the two parts are orthogonal in $L^2(\mu)$. Suppose that
$(z-\A)h$ is small compared with $h$ for some $z$ with $\Imm z\ne0$, so
that $h$ is almost an eigenfunction with the eigenvalue $z$. Then
\eqref{eq:imag-res} shows that $\nm{(\id-\J)h}^2$ is close to $2\nm h^2$,
which means that the odd part of $h$ carries about one half of $\nm h^2$.
By \eqref{eq:odd} the dissipation $-\Ree\ip{\A h}h$ is then at least about
$\tfrac\gamma2\nm h^2$, whereas it is close to $-\Ree z\,\nm h^2$ because
$\A h$ is close to $zh$. This is possible only if $\Ree z\le-\gamma/2$, up
to a small error. For exact eigenfunctions there is no error, and we obtain
the following statement.

\begin{corollary}\label{cor:nonreal}
If $\A h=\lambda h$ with $h\ne0$ and $\Imm\lambda\ne0$, then
$\Ree\lambda\le-\gamma/2$.
\end{corollary}
\begin{proof}
Equation~\eqref{eq:Jreal} gives $\Imm\lambda\,\ip h{\J h}=0$, hence
$\ip h{\J h}=0$ and, by \eqref{eq:Jform},
$\nm{(\id-\J)h}^2=2\nm h^2$. Inserting this into \eqref{eq:odd} yields
$-\Ree\lambda\,\nm h^2\ge\tfrac\gamma2\nm h^2$.
\end{proof}

\begin{remark}
The a priori argument above works in the whole half-plane
$\Ree z>-\gamma/2$, as \eqref{eq:apriori} shows; the smaller strip
$\Ree z\ge-\gamma/4$ was chosen only because it gives a uniform constant.
Nothing is asserted about the real axis: real spectrum may be present on
$(-\infty,0]$ and may accumulate at zero. It is the
resolvent norm near the imaginary axis, rather than the location of
eigenvalues alone, which controls the semigroup estimates of
Proposition~\ref{prop:reg}.
\end{remark}

\begin{remark}
\label{rem:PTresolvent}
To the best of our knowledge, the resolvent bound of
Proposition~\ref{prop:res} is new. The resolvent estimates for
kinetic Fokker--Planck operators near the imaginary axis which we are aware
of \cite[Theorem~4.1]{HerauNier}, \cite[Chapters~5 and~6]{HelfferNier},
\cite[Theorem~1.15]{ABL}, \cite{Wang} are derived from hypoelliptic
estimates, require growth or decay assumptions on the potential, and give no
information when $|\Imm z|$ is small. Proposition~\ref{prop:res} holds for
every smooth potential, in a whole strip to the left of the imaginary axis
and down to the real axis, with a constant which depends only on $\gamma$
and on $|\Imm z|$. The proof uses no hypoelliptic estimate, and its two steps are elementary:
pairing the resolvent equation with $h$ controls the velocity-odd part of
$h$ through the Gaussian Poincar\'e inequality, and pairing it with $\J h$
controls $|\Imm z|\,\ip{\J h}h$ through the identity \eqref{eq:Jreal},
which is a consequence of the adjoint relation \eqref{eq:J}. 
\end{remark}

We now deduce from the resolvent bound of
Proposition~\ref{prop:res} a decay estimate for the operators $\A^mP_t$,
composed with a power of the resolvent $(\id-\A)^{-1}$.

\begin{proposition}[Regularized decay estimate]\label{prop:reg}
Let $L\in C^\infty(\R^d)$. For every integer $m\ge1$ and every $t>0$,
\begin{equation}\label{eq:reg}
 \opnm{\A^mP_t(\id-\A)^{-m-2}}
          \le N_m\,t^{-m},\qquad
 N_m=\frac{3\cdot2^{(m+1)/2}}{\pi}\,(m-1)!
     +\frac{24}{\gamma(4+\gamma)}\Big(\frac{4m}{e\gamma}\Big)^m .
\end{equation}
\end{proposition}

Like Proposition~\ref{prop:res}, this estimate holds for every smooth
potential. The passage from a resolvent bound to a decay estimate is known
at the abstract level: for a bounded semigroup on a Hilbert space whose
generator has the origin as its only spectral point on the imaginary axis,
Batty, Chill and Tomilov \cite[Theorem~7.6]{BCT} prove that a resolvent
bound of order $|s|^{-1}$ as $s\to0$ along the imaginary axis, together
with a uniform bound as $|s|\to\infty$, is equivalent to
$\opnm{P_t\A(\id-\A)^{-1}}=O(t^{-1})$, and hence, by iteration, implies
$\opnm{\A^mP_t(\id-\A)^{-m}}=O(t^{-m})$. What is new in
Proposition~\ref{prop:reg} is that the resolvent bound~\eqref{eq:Rbound} is available for
every smooth potential, and that the constant $N_m$ is explicit, whereas
the theorem of \cite{BCT} gives no value of the constant. The explicit
constant comes from the direct proof which we now give. This proof
requires two additional powers of the resolvent, that is, the factor
$(\id-\A)^{-2}$ in \eqref{eq:reg}; they are harmless, since every power of
the resolvent is removed in Section~\ref{sec:proofmain}.

\begin{proof}
Let $\Gamma$ be the contour which travels upwards along the vertical line
$\Ree z=-\gamma/4$ from $-\gamma/4-i\infty$ to $-\gamma/4-i\gamma/4$, then
along the ray $z=-r-ir$ with $r$ decreasing from $\gamma/4$ to $0$, then
along the ray $z=-r+ir$ with $r$ increasing from $0$ to $\gamma/4$, and
finally upwards along $\Ree z=-\gamma/4$ from $-\gamma/4+i\gamma/4$ to
$-\gamma/4+i\infty$; see Figure~\ref{fig:contour}. Every point of $\Gamma$
other than the origin has
nonzero imaginary part and real part at least $-\gamma/4$, and therefore
lies in the resolvent set by Proposition~\ref{prop:res}. We claim that, for
$t>0$,
\begin{equation}\label{eq:contour}
 \A^mP_t(\id-\A)^{-m-2}
 =\frac1{2\pi i}\int_{\Gamma}
 e^{tz}\frac{z^m}{(1-z)^{m+2}}(z-\A)^{-1}\dd z ,
\end{equation}
the integral being absolutely convergent in operator norm, that is,
\[
 \int_{\Gamma}|e^{tz}|\,\frac{|z|^m}{|1-z|^{m+2}}\,
        \opnm{(z-\A)^{-1}}\,|\dd z|<\infty .
\]

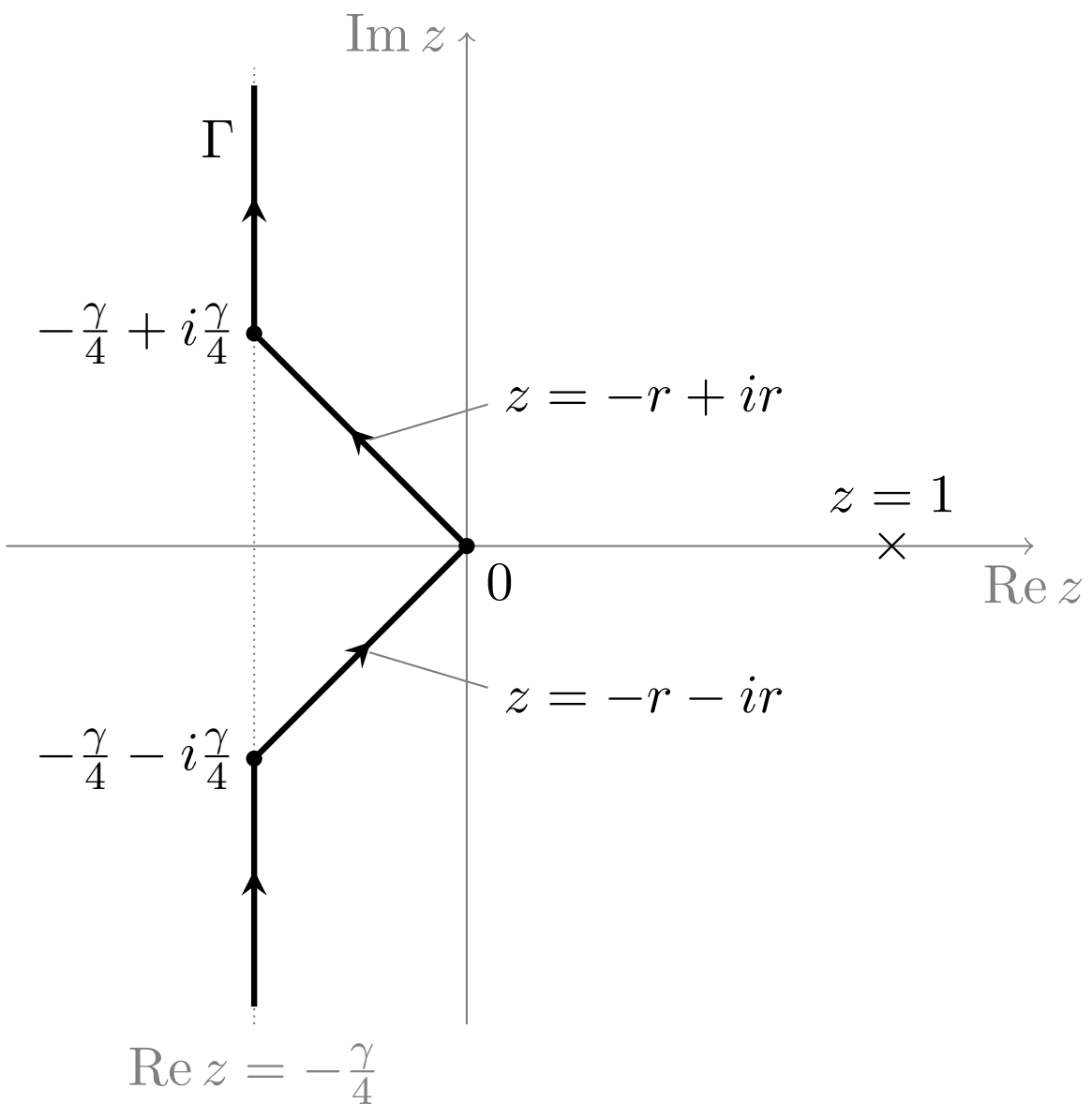
\begin{figure}[ht]
\centering
\scalebox{\fitfigscale}{%
\begin{tikzpicture}[scale=1.25,
  arr/.style={postaction={decorate},
     decoration={markings,mark=at position 0.55 with {\arrow{stealth}}}}]
 \draw[->,gray] (-2.6,0)--(3.2,0) node[below] {$\Ree z$};
 \draw[->,gray] (0,-2.7)--(0,2.9) node[left] {$\Imm z$};
 \draw[gray,densely dotted] (-1.2,-2.7)--(-1.2,2.7);
 \draw[very thick,arr] (-1.2,-2.6)--(-1.2,-1.2);
 \draw[very thick,arr] (-1.2,-1.2)--(0,0);
 \draw[very thick,arr] (0,0)--(-1.2,1.2);
 \draw[very thick,arr] (-1.2,1.2)--(-1.2,2.6);
 \fill (0,0) circle (1.3pt) node[below right] {$0$};
 \fill (-1.2,1.2) circle (1.3pt) node[left] {$-\tfrac\gamma4+i\tfrac\gamma4$};
 \fill (-1.2,-1.2) circle (1.3pt) node[left] {$-\tfrac\gamma4-i\tfrac\gamma4$};
 \draw (2.4,0) node {$\times$} node[above=3pt] {$z=1$};
 \node[left] at (-1.2,2.3) {$\Gamma$};
 \node[right] at (0.1,0.85) {$z=-r+ir$};
 \draw[gray,thin] (0.12,0.8)--(-0.55,0.6);
 \node[right] at (0.1,-0.85) {$z=-r-ir$};
 \draw[gray,thin] (0.12,-0.8)--(-0.55,-0.6);
 \node[below,gray] at (-1.2,-2.7) {$\Ree z=-\tfrac\gamma4$};
\end{tikzpicture}}
\caption{The contour $\Gamma$ used in the proof of
Proposition~\ref{prop:reg}. It follows the line $\Ree z=-\gamma/4$, except
between the points $-\gamma/4\pm i\gamma/4$, where it passes through the
origin along the two rays $z=-r\pm ir$. The cross marks the pole $z=1$ of the
factor $(1-z)^{-m-2}$, which lies to the right of $\Gamma$.}
\label{fig:contour}
\end{figure}

\textbf{The representation on a vertical line.} Fix a number $c$ with
$0<c<1$. Every point $z$ of the vertical line $\Ree z=c$ belongs to
$\rho(\A)$, and $\opnm{(z-\A)^{-1}}\le1/c$ there, by \eqref{eq:halfplane}.
Our starting point is a classical formula which expresses a semigroup
through the resolvent of its generator
\cite[Corollary~III.5.15]{EngelNagel}: for every $g\in\Dom(\A)$ and every
$t>0$,
\[
 P_tg=\frac1{2\pi i}\lim_{R\to\infty}\int_{c-iR}^{c+iR}
     e^{tz}(z-\A)^{-1}g\dd z .
\]
The integral is taken along the segment of the line $\Ree z=c$ between the
points $c-iR$ and $c+iR$. In \cite{EngelNagel} the formula is stated for a
general semigroup, under the condition that $\opnm{P_t}\le Me^{\omega t}$
for all $t\ge0$, with some constants $M\ge1$ and $\omega<c$. Here
$\opnm{P_t}\le1$, so that this condition holds with $M=1$ and $\omega=0$,
and every $c>0$ is admissible. Now let $j\ge1$ be an integer and
let $f\in L^2(\mu)$. The function $g=(\id-\A)^{-j}f$ belongs to $\Dom(\A)$,
because $(\id-\A)^{-1}$ maps $L^2(\mu)$ into $\Dom(\A)\subset L^2(\mu)$, and the formula
gives
\[
 P_t(\id-\A)^{-j}f=\frac1{2\pi i}\lim_{R\to\infty}\int_{c-iR}^{c+iR}
     e^{tz}(z-\A)^{-1}(\id-\A)^{-j}f\dd z,\qquad t>0.
\]
We now bring this formula into the form \eqref{eq:contour}, in three steps.

\emph{Step 1: reduction to the operators $P_t(\id-\A)^{-j}$.} Writing
$\A=\id-(\id-\A)$ and expanding the $m$th power by the binomial formula, we
obtain $\A^m=\sum_{l=0}^m(-1)^l\binom ml(\id-\A)^l$, and therefore
\[
 \A^m(\id-\A)^{-m-2}
       =\sum_{l=0}^m(-1)^l\binom ml(\id-\A)^{-(m+2-l)} .
\]
Consequently
\begin{equation}\label{eq:binomialsum}
 \A^mP_t(\id-\A)^{-m-2}f
   =\sum_{l=0}^m(-1)^l\binom ml\,P_t(\id-\A)^{-(m+2-l)}f ,
\end{equation}
and each term on the right-hand side is of the form $P_t(\id-\A)^{-j}f$,
with $j=m+2-l\ge2$, to which the formula above applies.

\emph{Step 2: the resolvent alone under the integral.} In the formula for
$P_t(\id-\A)^{-j}f$ the integrand contains the product
$(z-\A)^{-1}(\id-\A)^{-j}$. We replace it by means of the identity
\begin{equation}\label{eq:partialfractions}
 (z-\A)^{-1}(\id-\A)^{-j}
 =\frac{(z-\A)^{-1}}{(1-z)^j}
       -\sum_{k=1}^j\frac{(\id-\A)^{-k}}{(1-z)^{j+1-k}},
 \qquad j\ge1,
\end{equation}
valid for $z\in\rho(\A)$ with $z\ne1$. It is proved by induction on $j$.
For $j=1$ it reads
\[
 (z-\A)^{-1}(\id-\A)^{-1}
   =\frac{(z-\A)^{-1}-(\id-\A)^{-1}}{1-z},
\]
which is the resolvent identity; to check it, one multiplies both sides by
$(z-\A)$ on the left and by $(\id-\A)$ on the right, and uses
$(\id-\A)-(z-\A)=1-z$. If the identity holds for some $j\ge1$, we multiply
it by $(\id-\A)^{-1}$ and replace the product
$(z-\A)^{-1}(\id-\A)^{-1}$ in the first term by the right-hand side of the
case $j=1$; this gives the identity for $j+1$.

When \eqref{eq:partialfractions} is inserted in the formula for
$P_t(\id-\A)^{-j}f$, the terms of the sum over $k$ contribute nothing.
Indeed, in each of them the operator $(\id-\A)^{-k}$ does not depend on $z$,
so that the term is the fixed vector $(\id-\A)^{-k}f$ multiplied by the
number
\[
 \frac1{2\pi i}\lim_{R\to\infty}\int_{c-iR}^{c+iR}
        \frac{e^{tz}}{(1-z)^{n}}\dd z ,\qquad n=j+1-k\ge1,
\]
and this number is zero for $t>0$. This follows from the following
result on integrals of rational functions
\cite[Chapter~4, Section~5.3, case~3]{Ahlfors}: if $Q$ is a rational
function which has no pole on the real axis and tends to zero at infinity,
then
\[
 \lim_{X\to\infty}\int_{-X}^{X}Q(x)e^{ix}\dd x
   =2\pi i\sum\operatorname{Res}\big(Q(w)e^{iw}\big),
\]
where the sum is taken over the poles of $Q$ in the upper half-plane
$\Imm w>0$. To apply it we parametrize the line $\Ree z=c$ by $z=c+ix/t$
with $x$ real, so that $e^{tz}=e^{tc}e^{ix}$ and $\dd z=(i/t)\dd x$, and we
obtain
\[
 \frac1{2\pi i}\int_{c-iR}^{c+iR}\frac{e^{tz}}{(1-z)^{n}}\dd z
  =\frac{e^{tc}}{2\pi t}\int_{-tR}^{tR}Q(x)e^{ix}\dd x,
 \qquad Q(x)=\Big(1-c-\frac{ix}{t}\Big)^{-n}.
\]
The rational function $Q$ tends to zero at infinity, and its only pole is
the point $w=-i\,t(1-c)$, which lies in the lower half-plane because $t>0$
and $c<1$. There is therefore no pole in the upper half-plane, the sum of
residues is empty, and the limit as $R\to\infty$ is zero. We have thus
shown that
\[
 P_t(\id-\A)^{-j}f=\frac1{2\pi i}\lim_{R\to\infty}\int_{c-iR}^{c+iR}
     \frac{e^{tz}}{(1-z)^j}\,(z-\A)^{-1}f\dd z,\qquad t>0,\quad j\ge1.
\]

\emph{Step 3: summation.} We insert the last formula, with $j=m+2-l$, in
each term of \eqref{eq:binomialsum}. The scalar factors add up to
\[
 \sum_{l=0}^m(-1)^l\binom ml\frac1{(1-z)^{m+2-l}}
   =\frac1{(1-z)^{m+2}}\sum_{l=0}^m\binom ml\big(-(1-z)\big)^l
   =\frac{\big(1-(1-z)\big)^m}{(1-z)^{m+2}}
   =\frac{z^m}{(1-z)^{m+2}},
\]
again by the binomial formula. This gives \eqref{eq:contour} with the line
$\Ree z=c$ in place of $\Gamma$. On this line the resolvent is bounded by
$1/c$ and $|z^m(1-z)^{-m-2}|\le C|z|^{-2}$ for large $|z|$, so that the
integrand is $O(|\Imm z|^{-2})$ in operator norm and the integral converges
absolutely; the limit $R\to\infty$ is therefore an ordinary integral over
the whole line.

\textbf{Deformation of the contour.} The representation on the line
$\Ree z=c$ is of no use for large times, because $|e^{tz}|=e^{tc}$ on this
line, which grows exponentially with $t$. To obtain decay we must move the
path of integration into the half-plane $\Ree z<0$, where
$|e^{tz}|=e^{t\Ree z}$ is small. The path cannot be moved across the
negative real axis, which may contain spectrum up to the origin, and this
is the reason for the shape of $\Gamma$: far from the real axis it follows
the line $\Ree z=-\gamma/4$, and near the real axis it is forced to pass
through the origin, which it does along the two rays $z=-r\pm ir$. The
region between the line $\Ree z=c$
and $\Gamma$ consists of points which either have positive real part, or
have real part in $[-\gamma/4,0]$ and nonzero imaginary part; it is
therefore contained in $\rho(\A)$, and it does not contain the pole $z=1$
of the scalar factor. To deform the line into $\Gamma$ we first truncate at
height $\pm R$ and replace the corner of $\Gamma$ at the origin by the arc
of the circle $|z|=\varepsilon$ which joins the two rays through the right
half-plane; see Figure~\ref{fig:deformation}. Together with the two
horizontal segments at height $\pm R$, the truncated line and the modified
contour form a closed path, and the integrand of \eqref{eq:contour} is
holomorphic on the region enclosed by this path and on the path itself:
the resolvent $(z-\A)^{-1}$ is a holomorphic function of $z$ on
$\rho(\A)$, and the scalar factor is holomorphic away from $z=1$. By
Cauchy's theorem the integral along the closed path is zero. It remains to let $R\to\infty$ and
$\varepsilon\to0$, and to show that the contributions of the horizontal
segments and of the arc tend to zero. On the
horizontal segments at height $\pm R$, which have the fixed length
$c+\gamma/4$, the resolvent is bounded by \eqref{eq:Rbound}, which
applies because these segments lie in the region $\Ree z\ge-\gamma/4$,
$\Imm z\ne0$, and the scalar factor is $O(R^{-2})$, so their
contributions vanish as $R\to\infty$. On the arc $|z|=\varepsilon$ one has
$\opnm{(z-\A)^{-1}}=O(\varepsilon^{-1})$: where
$\Ree z\ge\varepsilon/\sqrt2$ the bound \eqref{eq:halfplane} by
$1/\Ree z$ suffices, and on the remaining part of the arc
$|\Imm z|\ge\varepsilon/\sqrt2$, so that \eqref{eq:Rbound} applies. Since
$|z^m|=\varepsilon^m$ and the arc has length at most $2\pi\varepsilon$, the
arc integral is $O(\varepsilon^m)$ for fixed $t>0$ and vanishes as
$\varepsilon\to0$, because $m\ge1$. The estimates below show that the
integral over $\Gamma$ converges absolutely, which completes the proof of
\eqref{eq:contour}.

\fitfigure{%
\begin{tikzpicture}[x=1.15cm,y=1.15cm,
  arr/.style={postaction={decorate},
     decoration={markings,mark=at position 0.55 with {\arrow{stealth}}}}]
 \fill[gray!15] (-1.3,-2.7)--(-1.3,-1.3)--(-0.283,-0.283)
      arc[start angle=-135,end angle=135,radius=0.4]
      --(-1.3,1.3)--(-1.3,2.7)--(1.7,2.7)--(1.7,-2.7)--cycle;
 \draw[->,gray] (-2.5,0)--(3.9,0) node[below] {$\Ree z$};
 \draw[->,gray] (0,-3.3)--(0,3.4) node[left] {$\Imm z$};
 \draw[thick,dashed,arr] (1.7,-2.7)--(1.7,2.7);
 \draw[very thick,arr] (-1.3,-2.7)--(-1.3,-1.3);
 \draw[very thick,arr] (-1.3,-1.3)--(-0.283,-0.283);
 \draw[very thick] (-0.283,-0.283)
      arc[start angle=-135,end angle=135,radius=0.4];
 \draw[very thick,arr] (-0.283,0.283)--(-1.3,1.3);
 \draw[very thick,arr] (-1.3,1.3)--(-1.3,2.7);
 \draw[thick] (-1.3,2.7)--(1.7,2.7);
 \draw[thick] (-1.3,-2.7)--(1.7,-2.7);
 \node[above] at (0.85,2.7) {$\Imm z=R$};
 \node[below] at (0.85,-2.7) {$\Imm z=-R$};
 \node[right] at (1.7,1.6) {$\Ree z=c$};
 \node[left] at (-1.3,2.1) {$\Gamma$};
 \node[right] at (0.5,-0.95) {$|z|=\varepsilon$};
 \draw[gray,thin] (0.62,-0.8)--(0.33,-0.23);
 \draw (3.1,0) node {$\times$} node[above=3pt] {$z=1$};
 \node[below left,gray] at (-1.3,0) {$-\tfrac\gamma4$};
 \fill[gray] (-1.3,0) circle (1pt);
 \draw[->,thick,gray!70!black] (1.0,1.25) arc[start angle=-60,end angle=240,radius=0.3];
\end{tikzpicture}}{%
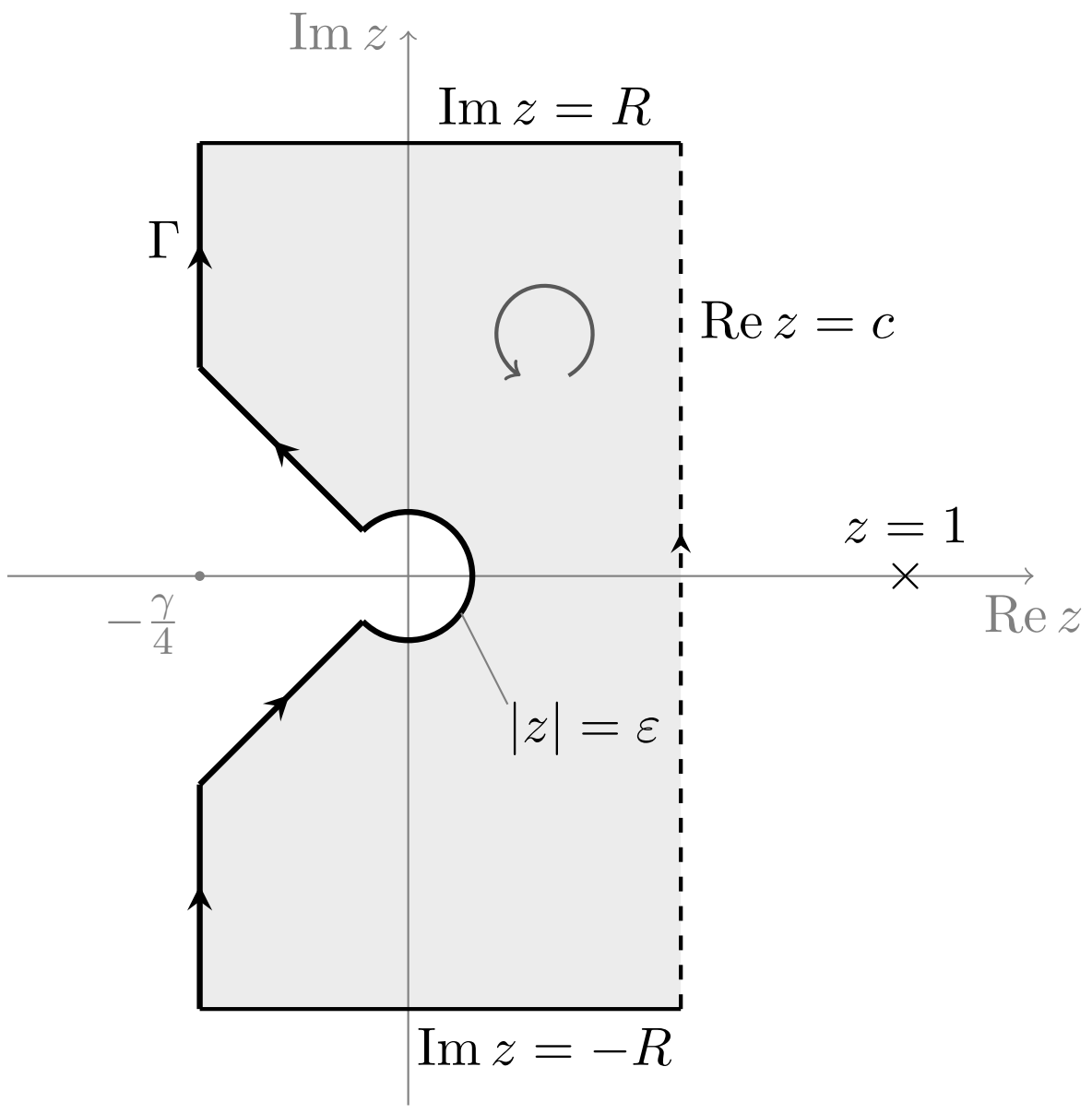
\caption{Deformation of the contour. The line $\Ree z=c$ (dashed) and the
contour $\Gamma$ are truncated at the heights $\Imm z=\pm R$, and the corner
of $\Gamma$ at the origin is replaced by an arc of the circle
$|z|=\varepsilon$. The shaded region between the two paths is contained in
$\rho(\A)$ and does not contain the pole $z=1$. The circular arrow indicates
the counterclockwise direction in which the boundary of this region is
described when Cauchy's theorem is applied.}
\label{fig:deformation}}

\textbf{Estimate of the integral.} On the two low-frequency rays
$z=-r\pm ir$, $0<r\le\gamma/4$, we have $|z|=\sqrt2\,r$, $|1-z|\ge1$,
$|e^{tz}|=e^{-tr}$ and, by \eqref{eq:Rbound},
\[
 \opnm{(z-\A)^{-1}}\le\frac4\gamma+\frac2r\le\frac3r .
\]
Since $|\dd z|=\sqrt2\dd r$ on each ray, the total contribution of the two
rays to \eqref{eq:contour} is therefore bounded in norm by
\[
 \frac{2}{2\pi}\int_0^{\gamma/4}e^{-tr}\,(\sqrt2\,r)^m\,\frac3r\,\sqrt2\dd r
 =\frac{3\cdot2^{(m+1)/2}}{\pi}\int_0^{\gamma/4} e^{-tr}r^{m-1}\dd r
       \le\frac{3\cdot2^{(m+1)/2}}{\pi}\,(m-1)!\,t^{-m}.
\]
On the two vertical parts $z=-\gamma/4+i\beta$, $|\beta|\ge\gamma/4$, we use
$|z|\le|1-z|$ to bound
\[
 \frac{|z|^m}{|1-z|^{m+2}}
       \le\frac1{(1+\gamma/4)^2+\beta^2},\qquad
 \opnm{(z-\A)^{-1}}\le\frac4\gamma+\frac2{|\beta|}\le\frac{12}\gamma ,
\]
and $|e^{tz}|=e^{-\gamma t/4}$, so that these parts contribute at most
\[
 \frac{2}{2\pi}\,\frac{12}\gamma\,e^{-\gamma t/4}
 \int_0^\infty\frac{\dd\beta}{(1+\gamma/4)^2+\beta^2}
 =\frac{24}{\gamma(4+\gamma)}\,e^{-\gamma t/4}
 \le\frac{24}{\gamma(4+\gamma)}\Big(\frac{4m}{e\gamma}\Big)^m t^{-m},
\]
where the last inequality uses that $t^me^{-\gamma t/4}$ attains its
maximum at $t=4m/\gamma$. Adding the two contributions, we obtain
\eqref{eq:reg} for every $t>0$, with the constant $N_m$ stated there.
\end{proof}

\begin{remark}
\label{rem:regularization}
The factor $(\id-\A)^{-m-2}$ in \eqref{eq:reg} cannot be dispensed with at
this stage. The exponent $m+2$ is the smallest one
for which the proof of Proposition~\ref{prop:reg} applies: on the vertical
part of the contour the integrand decays only through the scalar factor
$z^m(1-z)^{-m-2}$, which is of order $|z|^{-2}$ and hence integrable,
whereas with $m$ or $m+1$ powers of the resolvent it would be of order $1$
or $|z|^{-1}$. The factor is removed in
Section~\ref{sec:proofmain} by means of Proposition~\ref{prop:smooth},
which is a regularity property of the semigroup at positive times and is
the place where Assumption~\ref{ass:WH} is used.
\end{remark}

\begin{remark}\label{rem:BCT}
Wang \cite{Wang} and Pan, Wang and Zhu \cite{PWZ} also write the semigroup
as an integral of the resolvent along a contour, but with a different aim.
They study the solution $P_th$ itself, for potentials whose gradient
decays at infinity, and they compute its leading term as $t\to\infty$.
\end{remark}

\subsection{Smoothing at positive times}\label{sec:energy}

This subsection uses Assumption~\ref{ass:WH} to remove the regularization
$(\id-\A)^{-m-2}$ in Proposition~\ref{prop:reg}.

\begin{lemma}[A bound for the generator]\label{lem:mult}
If $h\in L^2(\mu)$ and the right-hand side of \eqref{eq:Agraph}
is finite, then $h\in\Dom(\A)$ and
\begin{equation}\label{eq:Agraph}
 \nm{\A h}^2\le C_{\rm mult}\big(
 \nm{\nabla_vh}^2+\nm{\nabla_xh}^2+
 \nm{\nabla_v^2h}^2+\nm{\nabla_x\nabla_vh}^2\big),
\end{equation}
with a constant $C_{\rm mult}$ depending only on $\gamma,\sigma$ and
$C_L$, whose value is given in Table~\ref{tab:constants}.
\end{lemma}
\begin{proof}
Let first $h\in\Cc$. We use the representation \eqref{eq:adjointform} of
$\Lham h$ and $\LFD h$ in terms of the formal adjoints $\nabla_v^*$ and
$\nabla_x^*$ defined in \eqref{eq:adjoints}.

Let $F=(F_1,\dots,F_d)$ be a vector field with components in $\Cc$. By
\eqref{eq:adjoints} we have $\nabla_v^*F=\sum_i\partial_{v_i}^*F_i$ with
$\partial_{v_i}^*=-\partial_{v_i}+v_i/\sigma$, and these operators satisfy
$[\partial_{v_j},\partial_{v_i}^*]=\delta_{ij}/\sigma$. Hence
\[
 \nm{\nabla_v^*F}^2
 =\sum_{i,j}\ip{\partial_{v_j}\partial_{v_i}^*F_i}{F_j}
 =\sum_{i,j}\ip{\partial_{v_j}F_i}{\partial_{v_i}F_j}
        +\frac1\sigma\nm F^2
 \le\nm{\nabla_vF}^2+\frac1\sigma\nm F^2 ,
\]
where the last step is the Cauchy--Schwarz inequality. In the same way,
$\nabla_x^*F=\sum_i\partial_{x_i}^*F_i$ with
$\partial_{x_i}^*=-\partial_{x_i}+\partial_{x_i}L/\sigma$, and since
$[\partial_{x_j},\partial_{x_i}^*]
=\partial_{x_i}\partial_{x_j}L/\sigma$,
\[
 \nm{\nabla_x^*F}^2
 =\sum_{i,j}\ip{\partial_{x_j}F_i}{\partial_{x_i}F_j}
        +\frac1\sigma\ip{(\nabla^2L)F}{F}
 \le\nm{\nabla_xF}^2+\frac1\sigma\nm{(\nabla^2L)F}\,\nm F .
\]
We apply the first inequality to $F=\nabla_xh$ and to $F=\nabla_vh$, and
the second one to $F=\nabla_vh$. In the second one we use
Assumption~\ref{ass:WH}, applied to the components of $\nabla_vh$, which
gives $\nm{(\nabla^2L)\nabla_vh}^2\le
C_L(\nm{\nabla_vh}^2+\nm{\nabla_x\nabla_vh}^2)$, and the elementary
inequality $(A^2+B^2)^{1/2}A\le\tfrac32A^2+\tfrac12B^2$. This gives
\begin{align*}
 \nm{\nabla_v^*\nabla_xh}^2
 &\le\nm{\nabla_v\nabla_xh}^2+\frac1\sigma\nm{\nabla_xh}^2,\\
 \nm{\nabla_v^*\nabla_vh}^2
 &\le\nm{\nabla_v^2h}^2+\frac1\sigma\nm{\nabla_vh}^2,\\
 \nm{\nabla_x^*\nabla_vh}^2
 &\le\Big(1+\frac{\sqrt{C_L}}{2\sigma}\Big)\nm{\nabla_x\nabla_vh}^2
     +\frac{3\sqrt{C_L}}{2\sigma}\nm{\nabla_vh}^2 .
\end{align*}
By \eqref{eq:adjointform}, $\A h=-\Lham h+\gamma\LFD h$ is the sum of these
three expressions multiplied by $-\sigma$, $-\gamma\sigma$ and $\sigma$
respectively. Since $\nm{f_1+f_2+f_3}^2\le3\sum_i\nm{f_i}^2$, we conclude
that
\[
 \nm{\A h}^2\le3\sigma\nm{\nabla_xh}^2
 +3\sigma\Big(\frac{3\sqrt{C_L}}2+\gamma^2\Big)\nm{\nabla_vh}^2
 +3\gamma^2\sigma^2\nm{\nabla_v^2h}^2
 +3\Big(2\sigma^2+\frac{\sigma\sqrt{C_L}}2\Big)\nm{\nabla_x\nabla_vh}^2 ,
\]
which is \eqref{eq:Agraph} for test functions, with
\[
 C_{\rm mult}=3\max\Big\{\sigma\Big(\frac{3\sqrt{C_L}}2+\gamma^2\Big),\
 \gamma^2\sigma^2,\ 2\sigma^2+\frac{\sigma\sqrt{C_L}}2\Big\};
\]
the coefficient $3\sigma$ of $\nm{\nabla_xh}^2$ does not exceed the first
term of the maximum because $C_L\ge1$.

Let now $h\in L^2(\mu)$ be such that the right-hand side of
\eqref{eq:Agraph} is finite. By truncation and mollification there are
functions $h_n\in\Cc$ which converge to $h$ in $L^2(\mu)$ and whose
derivatives $\nabla_vh_n$, $\nabla_xh_n$, $\nabla_v^2h_n$ and
$\nabla_x\nabla_vh_n$ converge in $L^2(\mu)$ to the corresponding
derivatives of $h$. The bound just proved, applied to the differences
$h_n-h_k$, shows that $\A h_n$ is a Cauchy sequence in $L^2(\mu)$. Since
$\A$ is a closed operator, $h$ belongs to $\Dom(\A)$ and $\A h$ is the
limit of $\A h_n$, and \eqref{eq:Agraph} for $h$ follows by passing to the
limit.
\end{proof}

The two energy inequalities involve three positive coefficients $a,b,c$,
which depend only on $\gamma$, $\sigma$ and $C_L$, satisfy $ac=4b^2$, and
are given in Table~\ref{tab:constants}. For $w\in H^1(\mu)$ we define the
energy
\begin{equation}\label{eq:staticenergy}
 \cG(w)=\nm{w}^2+a\nm{\nabla_vw}^2
 +2b\Ree\ip{\nabla_vw}{\nabla_xw}+c\nm{\nabla_xw}^2,
\end{equation}
and, for $w\in H^2(\mu)$, the dissipation
\begin{equation}\label{eq:firstdiss}
 \cD(w)=\nm{\nabla_vw}^2+\nm{\nabla_xw}^2+
 \nm{\nabla_v^2w}^2+\nm{\nabla_x\nabla_vw}^2 .
\end{equation}
Since $2b=\sqrt{ac}$, the mixed term of $\cG$ is bounded by one half of the
sum of the two gradient terms, so that $\cG(w)$ is equivalent to the
squared norm of $w$ in $H^1(\mu)$:
\begin{equation}\label{eq:Gequiv}
 \tfrac12\min(1,a,c)\big(\nm w^2+\nm{\nabla_{x,v}w}^2\big)\le\cG(w)\le
 \tfrac32\max(1,a,c)\big(\nm w^2+\nm{\nabla_{x,v}w}^2\big).
\end{equation}

\begin{lemma}[Regularization with norm loss]\label{lem:reverse}
For $h\in L^2(\mu)$ and $0<t\le1$ we have $P_th\in H^1(\mu)$ and
\begin{equation}\label{eq:reverse-t}
 \frac a2\,t\,\nm{\nabla_vP_th}^2+
 \frac c2\,t^3\nm{\nabla_xP_th}^2
       \le\nm h^2-\nm{P_th}^2.
\end{equation}
In particular,
\begin{equation}\label{eq:reverse-one}
 \nm{\nabla_{x,v}P_1h}^2
 \le2\max(a^{-1},c^{-1})\big(\nm h^2-\nm{P_1h}^2\big).
\end{equation}
\end{lemma}

\begin{lemma}[Static energy inequality]\label{lem:static}
For $h\in H^1(\mu)$ and $t\ge0$ we have $P_th\in H^1(\mu)$ and
\begin{equation}\label{eq:static}
 \cG(P_th)+c_*\int_0^t\cD(P_sh)\dd s\le\cG(h),
\end{equation}
with the constant $c_*$ of Table~\ref{tab:constants}.
In particular, $\nabla_v^2P_sh$ and $\nabla_x\nabla_vP_sh$ belong to
$L^2(\mu)$ for almost every $s>0$.
\end{lemma}

Lemmas~\ref{lem:reverse} and~\ref{lem:static} are classical results, due
to H\'erau \cite{Herau}, Villani \cite{Villani} and Zhang \cite{Zhang};
Appendix~\ref{apx:energy} provides a rigorous verification under our
assumptions, and the precise references are given there.

Lemma~\ref{lem:reverse} is used in Section~\ref{sec:proofmain} as follows:
applied to the function $\A^mP_{t-1}h_0$, the inequality
\eqref{eq:reverse-one} bounds its gradient at time $t$ by the decrease of
its norm between the times $t-1$ and $t$.

The next proposition states that $\A^jP_\tau$ is a bounded operator on
$L^2(\mu)$ for every $\tau>0$ and every $j\ge1$.

\begin{proposition}[Positive-time smoothing]\label{prop:smooth}
For $\tau>0$ and every integer $j\ge1$,
\begin{equation}\label{eq:fixedsmooth}
 P_\tau L^2(\mu)\subset\Dom(\A^j),\qquad
                         \opnm{\A^jP_\tau}<\infty.
\end{equation}
In particular,
\begin{equation}\label{eq:remove}
 \opnm{(\id-\A)^rP_1}<\infty,\qquad r\ge1.
\end{equation}
\end{proposition}
\begin{proof}
Let $0<\tau\le1$ and $h\in L^2(\mu)$. By Lemma~\ref{lem:reverse} at the
time $\tau/2$, the function $w=P_{\tau/2}h$ belongs to $H^1(\mu)$ and
\[
 \frac{a\tau}{4}\nm{\nabla_vw}^2+\frac{c\tau^3}{16}\nm{\nabla_xw}^2
 \le\nm h^2,
 \qquad\text{so that}\qquad
 a\nm{\nabla_vw}^2+c\nm{\nabla_xw}^2\le\frac{16}{\tau^3}\nm h^2
\]
because $\tau\le1$. Since the mixed term of $\cG$ is bounded by one half
of the sum of the two gradient terms, and since $\nm w\le\nm h$, this
gives
\[
 \cG(P_{\tau/2}h)\le\nm w^2
 +\tfrac32\big(a\nm{\nabla_vw}^2+c\nm{\nabla_xw}^2\big)
 \le\frac{25}{\tau^3}\nm h^2 .
\]
We may therefore
apply Lemma~\ref{lem:static} with the initial datum $P_{\tau/2}h$ on the
time interval $[0,\tau/2]$, which gives
\[
 c_*\int_{\tau/2}^\tau\cD(P_sh)\dd s\le\cG(P_{\tau/2}h)
 \le\frac{25}{\tau^3}\nm h^2 .
\]
For almost every $s\in(\tau/2,\tau)$ the four derivatives in
$\cD(P_sh)$ are thus in $L^2(\mu)$, so that Lemma~\ref{lem:mult} shows that
$P_sh\in\Dom(\A)$ and $\nm{\A P_sh}^2\le C_{\rm mult}\cD(P_sh)$.
Consequently
\[
 \int_{\tau/2}^\tau\nm{\A P_sh}^2\dd s
 \le\frac{C_{\rm mult}}{c_*}\,\cG(P_{\tau/2}h)
 \le\frac{25\,C_{\rm mult}}{c_*\tau^3}\nm h^2.
\]
If $h\in\Dom(\A)$, then $\A P_sh=P_s\A h$ and the contraction
property of $P_s$ makes $s\mapsto\nm{\A P_sh}$ nonincreasing, whence
\[
 \frac\tau2\,\nm{\A P_\tau h}^2\le\int_{\tau/2}^\tau\nm{\A P_sh}^2\dd s
   \le\frac{25\,C_{\rm mult}}{c_*\tau^3}\nm h^2 ,
 \qquad\text{that is,}\qquad
 \nm{\A P_\tau h}\le\frac K{\tau^2}\nm h
\]
with $K=(50\,C_{\rm mult}/c_*)^{1/2}$.
The density of $\Dom(\A)$ in $L^2(\mu)$ and the closedness of $\A$ extend
this bound to all $h\in L^2(\mu)$, and show in particular that
$P_\tau L^2(\mu)\subset\Dom(\A)$.
For $\tau>1$ we write $P_\tau h=P_1(P_{\tau-1}h)$, which belongs to
$\Dom(\A)$ by the case $\tau=1$, and since $P_{\tau-1}$ is a contraction
we get $\nm{\A P_\tau h}\le K\nm{P_{\tau-1}h}\le K\nm h$. Altogether,
\begin{equation}\label{eq:APtau}
 \nm{\A P_\tau h}\le K\max\Big(1,\frac1\tau\Big)^{2}\nm h,
 \qquad h\in L^2(\mu),\quad\tau>0,
\end{equation}
which proves \eqref{eq:fixedsmooth} for $j=1$.

Let now $j\ge2$, $\tau>0$ and
$s=\tau/j$. The operator $\A P_s$ is bounded, and it commutes with the
semigroup because $\A$ commutes with the semigroup $P_s$ on $\Dom(\A)$. For
$h\in L^2(\mu)$ we therefore have
\[
 \A P_\tau h=\A P_sP_{\tau-s}h=P_{\tau-s}(\A P_sh),
\]
and the right-hand side is again $P_{\tau-s}=P_{(j-1)s}$ applied to a
function of $L^2(\mu)$. Repeating this argument $j$ times, we obtain
$P_\tau h\in\Dom(\A^j)$ and
\begin{equation}\label{eq:powerinduction}
 \A^jP_\tau h=(\A P_{\tau/j})^jh,\qquad\text{hence}\qquad
 \nm{\A^jP_\tau h}\le\opnm{\A P_{\tau/j}}^j\,\nm h
 \le K^j\max\Big(1,\frac j\tau\Big)^{2j}\nm h
\end{equation}
by \eqref{eq:APtau} at the time $\tau/j$. This proves
\eqref{eq:fixedsmooth}.
Finally, the binomial expansion
\[
 (\id-\A)^rP_1=\sum_{j=0}^r\binom rj(-1)^j\A^jP_1
\]
proves \eqref{eq:remove}.
\end{proof}

Note that no bounds on higher derivatives of $L$ are required:
Assumption~\ref{ass:WH} suffices for all $j$.

\subsection{Conclusion of the proof}\label{sec:proofmain}

We combine Proposition~\ref{prop:reg} with the estimates of
Section~\ref{sec:energy}.

\textbf{Generator powers.} For $t\ge2$ we write
\begin{equation}\label{eq:remove-contour}
 \A P_t=
 \big[\A P_{t-1}(\id-\A)^{-3}\big]\,\big[(\id-\A)^3P_1\big],
\end{equation}
which is legitimate because $P_1L^2(\mu)\subset\Dom(\A^3)$ and the
resolvent commutes with the semigroup. By Proposition~\ref{prop:reg} with
$m=1$ and by \eqref{eq:remove} with $r=3$,
\[
 \opnm{\A P_t}\le N_1(t-1)^{-1}\opnm{(\id-\A)^3P_1}\le Ct^{-1},
 \qquad t\ge2 .
\]
For $1\le t\le2$ we use $\A P_t=(\A P_1)P_{t-1}$ and the contraction
property, so that altogether $\opnm{\A P_t}\le K_1t^{-1}$ for all $t\ge1$,
with the constant $K_1$ of Table~\ref{tab:constants}.

Let now $m\ge1$. For $t\ge m$ we have $\A^mP_t=(\A P_{t/m})^m$ by
\eqref{eq:powerinduction}, which is
bounded in norm by $(mK_1)^mt^{-m}$ since $t/m\ge1$. For $1\le t\le m$ we
use $\A^mP_t=(\A^mP_1)P_{t-1}$, which is bounded by $\opnm{\A^mP_1}$, and
we enlarge the constant by the factor $m^m$. With $K_0=1$ we have proved
\begin{equation}\label{eq:longnorm}
 \opnm{\A^mP_t}\le K_mt^{-m},\qquad m\ge0,\quad t\ge1,
\end{equation}
and squaring gives \eqref{eq:main-even} with $C_m=K_m^2$.

\textbf{Regularity at positive times.} The inclusion
$P_tL^2(\mu)\subset\Dom(\A^m)$ for $t>0$ is part of
Proposition~\ref{prop:smooth}. To see that $\A^mP_t$ maps $L^2(\mu)$ into
$H^1(\mu)$ for $t>0$, we set $s=\min\{1,t/2\}$ and write
$\A^mP_t=P_s\big(\A^mP_{t-s}\big)$; the inner operator is bounded on
$L^2(\mu)$ by Proposition~\ref{prop:smooth}, and $P_s$ maps $L^2(\mu)$
into $H^1(\mu)$ by Lemma~\ref{lem:reverse}, since $0<s\le1$.

\textbf{The full gradient.} For $m\ge0$ and $t\ge1$, the dissipation
identity \eqref{eq:dirichlet}, applied to $\A^mP_th_0\in\Dom(\A)$, gives
together with \eqref{eq:longnorm}
\begin{equation}\label{eq:velocitylong}
 \gamma\sigma\nm{\nabla_v\A^mP_th_0}^2
 =-\Ree\ip{\A^{m+1}P_th_0}{\A^mP_th_0}
 \le K_{m+1}K_mt^{-(2m+1)}\nm{h_0}^2.
\end{equation}
This identity alone does not control the position gradient. To recover it
we apply, for $t\ge2$, the inequality \eqref{eq:reverse-one} to the
function $\A^mP_{t-1}h_0$, and we write the loss of norm over the interval
$[t-1,t]$ as the integral of its derivative. Indeed, the function
$g=\A^mP_{t-1}h_0$ belongs to $\Dom(\A)$ by
Proposition~\ref{prop:smooth}, so that, by
Proposition~\ref{prop:maxdissipative}, the function
$s\mapsto P_{s-(t-1)}g=\A^mP_sh_0$ is continuously differentiable on
$[t-1,t]$ with values in $L^2(\mu)$, with derivative $\A^{m+1}P_sh_0$;
hence
\[
 \frac{\dd}{\dd s}\nm{\A^mP_sh_0}^2
 =2\Ree\ip{\A^{m+1}P_sh_0}{\A^mP_sh_0},
\]
and the fundamental theorem of calculus gives the second line of the
following estimate:
\begin{align*}
 \nm{\nabla_{x,v}\A^mP_th_0}^2
 &\le2\max(a^{-1},c^{-1})\big(\nm{\A^mP_{t-1}h_0}^2
                         -\nm{\A^mP_th_0}^2\big)\\
 &=-4\max(a^{-1},c^{-1})\Ree\int_{t-1}^t
            \ip{\A^{m+1}P_sh_0}{\A^mP_sh_0}\dd s\\
 &\le4\max(a^{-1},c^{-1})K_{m+1}K_m\nm{h_0}^2
                         \int_{t-1}^t s^{-(2m+1)}\dd s\\
 &\le2^{2m+3}\max(a^{-1},c^{-1})K_{m+1}K_m\,
                         t^{-(2m+1)}\nm{h_0}^2 ,
\end{align*}
where we used \eqref{eq:longnorm} for $s\ge t-1\ge1$ and $t-1\ge t/2$. This
proves \eqref{eq:main-odd} with
$D_m=2^{2m+3}\max(a^{-1},c^{-1})K_{m+1}K_m$ and completes the proof of
Theorem~\ref{thm:main}. \qed

\begin{table}[!htb]
\centering\small
\renewcommand{\arraystretch}{2.1}
\begin{tabular}{@{}lll@{}}
\toprule
Constant & Value & Where\\
\midrule
$C_L$ & constant of the relative bound on the Hessian
      & Ass.~\ref{ass:WH}\\
$N_m$ & $\dfrac{3\cdot2^{(m+1)/2}}{\pi}\,(m-1)!
        +\dfrac{24}{\gamma(4+\gamma)}\Big(\dfrac{4m}{e\gamma}\Big)^m$
      & Prop.~\ref{prop:reg}\\
$C_{\rm mult}$ & $3\max\big\{\sigma(\tfrac32\sqrt{C_L}+\gamma^2),\
 \gamma^2\sigma^2,\ 2\sigma^2+\tfrac12\sigma\sqrt{C_L}\big\}$
      & Lemma~\ref{lem:mult}\\
$c/b$ & $\min\Big\{\dfrac16,\ \dfrac{\gamma\sigma}{16(1+C_L)}\Big\}$
      & \eqref{eq:coeffexplicit}\\
$a/b$ & $4b/c$ & \eqref{eq:coeffexplicit}\\
$b$ & $\min\Big\{1,\ \gamma\sigma\Big[1+\dfrac ab+2\sqrt{C_L}
   +\dfrac{4C_L\,b}{\gamma\sigma\,c}
   +4\Big(\dfrac ab+2+\gamma+\dfrac cb\sqrt{C_L}\Big)^2\Big]^{-1}\Big\}$
      & \eqref{eq:coeffexplicit}\\
$c_*$ & $\min(\gamma\sigma,\,b,\,\gamma\sigma a,\,\gamma\sigma c/2)
        =\min(b,\,\gamma\sigma c/2)$
      & Lemma~\ref{lem:static}\\
$K$ & $(50\,C_{\rm mult}/c_*)^{1/2}$
      & Prop.~\ref{prop:smooth}\\
$K_1$ & $2\big(2+9\gamma^{-2}\big)(1+9K)^3$
      & \eqref{eq:longnorm}\\
$K_m$ & $(m^3K_1)^m$ for $m\ge1$, and $K_0=1$
      & \eqref{eq:longnorm}\\
$C_m$ & $K_m^2$
      & Thm.~\ref{thm:main}\\
$D_m$ & $2^{2m+3}\,c^{-1}K_{m+1}K_m$
      & Thm.~\ref{thm:main}\\
\bottomrule
\end{tabular}
\caption{The constants of the paper. They depend only on $m$,
$\gamma$, $\sigma$ and $C_L$.}
\label{tab:constants}
\end{table}

\begin{remark}\label{rem:constants}
The values of Table~\ref{tab:constants} are obtained by keeping track of
the constants in the proofs. They are upper bounds, and we do not claim
that they are optimal. The
coefficients $a,b,c$ satisfy $c\le a$ and $b\le\gamma\sigma$, which gives
the second expression of $c_*$ and $\max(a^{-1},c^{-1})=c^{-1}$. The value
of $C_{\rm mult}$ is obtained in the proof of Lemma~\ref{lem:mult}. The proof of
Proposition~\ref{prop:smooth} gives $\opnm{\A^jP_1}\le K^jj^{2j}$, hence
$\opnm{(\id-\A)^3P_1}\le1+3K+48K^2+729K^3\le(1+9K)^3$, and
Proposition~\ref{prop:reg} gives
$N_1=6/\pi+96/(e\gamma^2(4+\gamma))\le2+9\gamma^{-2}$. Inserting these
bounds into the argument of this subsection, and using
$(t-1)^{-1}\le2/t$ for $t\ge2$ and
$\opnm{\A^mP_1}\le K^mm^{2m}\le(m^3K)^mt^{-m}$ for $1\le t\le m$, we obtain
the values of $K_1$ and $K_m$.
\end{remark}

\section{Sharpness: the free dynamics}\label{sec:sharp}

This section gives the proof of Theorem~\ref{thm:sharp}. We set
$L=0$, so that
\[
 \dd\mu=e^{-|v|^2/(2\sigma)}\dd x\dd v,\qquad
 \A=-v\cdot\nabla_x+\gamma\sigma\Delta_v-\gamma v\cdot\nabla_v .
\]
The operator $\A$ is then the generator of the diffusion
\[
 \dd X_t=-V_t\dd t,\qquad
 \dd V_t=-\gamma V_t\dd t+\sqrt{2\gamma\sigma}\,\dd W_t ,
\]
in the sense that $(P_tf)(x,v)$ is the expectation of $f(X_t,V_t)$ when the
process starts from $(x,v)$. The first equation has the opposite sign
compared with the underdamped Langevin dynamics \eqref{eq:SDE}. The
reason is
that \eqref{eq:A} is the equation satisfied by the density of the law of
that dynamics, whereas a generator describes the expectations of functions
of the process; passing from one description to the other replaces the
operator by its adjoint, and this changes the sign of the transport term
$v\cdot\nabla_x$. Subtracting $\gamma^{-1}$ times the second equation
from the first one, we see that the drift cancels:
\[
 \dd\big(X_t-V_t/\gamma\big)=-\sqrt{2\sigma/\gamma}\,\dd W_t .
\]
Thus $X_t-V_t/\gamma$ is a Brownian motion with diffusivity
$\sigma/\gamma$, and if the process starts from $(x,v)$, the expectation
of $h(X_t-V_t/\gamma)$ is $(e^{(\sigma/\gamma)t\Delta}h)(x-v/\gamma)$. In
other words, the functions of $x-v/\gamma$ evolve by the heat equation.
To make this precise we define, for
$h\in L^2(\R^d,\dd x)$,
\[
 (\mathcal Ih)(x,v)=(2\pi\sigma)^{-d/4}\,h(x-v/\gamma).
\]
Integrating first in $x$, for fixed $v$, and then in $v$, we find
\[
 \nm{\mathcal Ih}^2
 =(2\pi\sigma)^{-d/2}\int_{\R^d}\Big(\int_{\R^d}|h(x-v/\gamma)|^2\dd x\Big)
       e^{-|v|^2/(2\sigma)}\dd v
 =\|h\|_{L^2(\dd x)}^2 ,
\]
because the inner integral does not depend on $v$ and the Gaussian
$e^{-|v|^2/(2\sigma)}$ has total mass $(2\pi\sigma)^{d/2}$. Thus
$\mathcal I$ is an isometry from $L^2(\R^d,\dd x)$ into $L^2(\mu)$, and its
range is a closed subspace of $L^2(\mu)$.

\begin{proposition}[An invariant heat subspace]\label{prop:heat}
Let $L=0$. For every $h\in L^2(\R^d,\dd x)$ and every $t\ge0$,
\begin{equation}\label{eq:heatsub}
 P_t\,\mathcal Ih=\mathcal I\,e^{(\sigma/\gamma)t\Delta}h .
\end{equation}
Thus the range of $\mathcal I$ is invariant under $P_t$, and on this
subspace $P_t$ is unitarily equivalent to the heat semigroup on $\R^d$
with diffusivity $\sigma/\gamma$.
\end{proposition}

\begin{proof}The proof which follows is a direct computation; it does not use the
stochastic representation of the semigroup.
Let first $h$ be a Schwartz function, and write $f(x,v)=h(x-v/\gamma)$.
Then $\nabla_xf=(\nabla h)(x-v/\gamma)$,
$\nabla_vf=-\gamma^{-1}(\nabla h)(x-v/\gamma)$ and
$\Delta_vf=\gamma^{-2}(\Delta h)(x-v/\gamma)$, so that the three terms of
$\A f$ are $-v\cdot\nabla_xf=-v\cdot\nabla h$,
$\gamma\sigma\Delta_vf=(\sigma/\gamma)\Delta h$ and
$-\gamma v\cdot\nabla_vf=v\cdot\nabla h$, all evaluated at the point
$x-v/\gamma$. The two transport terms cancel, and
$\A[h(x-v/\gamma)]=(\sigma/\gamma)(\Delta h)(x-v/\gamma)$.
The function $\mathcal Ih$ is smooth and belongs to $L^2(\mu)$, and so does
the right-hand side, since $\Delta h$ is again a Schwartz function. The
maximal-domain characterization \eqref{eq:maxdomain} therefore shows that
$\mathcal Ih\in\Dom(\A)$ and that
\begin{equation}\label{eq:intertwine}
 \A\,\mathcal Ih=\mathcal I\,\frac\sigma\gamma\Delta h
 \qquad\text{for every Schwartz function }h .
\end{equation}
Consider now the function $u_t=\mathcal Ie^{(\sigma/\gamma)t\Delta}h$,
$t\ge0$. For a Schwartz function $h$, the heat flow
$t\mapsto e^{(\sigma/\gamma)t\Delta}h$ consists of Schwartz functions and
is continuously differentiable with values in $L^2(\R^d,\dd x)$, with time
derivative $(\sigma/\gamma)\Delta e^{(\sigma/\gamma)t\Delta}h$. Since
$\mathcal I$ is an isometry, $t\mapsto u_t$ is continuously differentiable
with values in $L^2(\mu)$, and by \eqref{eq:intertwine} it takes its
values in $\Dom(\A)$ and satisfies
$\frac{\dd}{\dd t}u_t
=\mathcal I\,(\sigma/\gamma)\Delta e^{(\sigma/\gamma)t\Delta}h=\A u_t$,
with $u_0=\mathcal Ih$.
Thus $u_t$ is a solution of \eqref{eq:A} with initial datum
$\mathcal Ih\in\Dom(\A)$, and by the uniqueness statement of
Proposition~\ref{prop:maxdissipative} it coincides with
$P_t\mathcal Ih$. This proves \eqref{eq:heatsub} for
Schwartz functions. The general case follows because Schwartz functions
are dense in $L^2(\R^d,\dd x)$ and both sides of \eqref{eq:heatsub} depend
continuously on $h\in L^2(\R^d,\dd x)$: the semigroups $P_t$ and
$e^{(\sigma/\gamma)t\Delta}$ are contractions on $L^2(\mu)$ and on
$L^2(\R^d,\dd x)$ respectively, and $\mathcal I$ is an isometry.
\end{proof}

\begin{proof}[Proof of Theorem~\ref{thm:sharp}]
Let $h$ be a Schwartz function and let $t>0$ and $m\ge0$. Since the heat
flow and the Laplacian preserve the Schwartz class, we may apply
\eqref{eq:intertwine} $m$ times to the right-hand side of
\eqref{eq:heatsub}, which gives $P_t\mathcal Ih\in\Dom(\A^m)$ and
\begin{equation}\label{eq:heatpowers}
 \A^mP_t\,\mathcal Ih=\mathcal I\,g,\qquad
 g=\Big(\frac\sigma\gamma\Delta\Big)^me^{(\sigma/\gamma)t\Delta}h .
\end{equation}
We compute the norm of the right-hand side with the Fourier transform
$\widehat h(\xi)=\int e^{-ix\cdot\xi}h(x)\dd x$. The function $g$ has the
Fourier transform $(-s)^me^{-ts}\,\widehat h(\xi)$, where
$s=(\sigma/\gamma)|\xi|^2$, and by the Plancherel theorem and the isometry
property of $\mathcal I$,
\[
 \frac{\nm{\A^mP_t\,\mathcal Ih}^2}{\nm{\mathcal Ih}^2}
 =\frac{\|g\|_{L^2(\dd x)}^2}{\|h\|_{L^2(\dd x)}^2}
 =\frac{\int_{\R^d}s^{2m}e^{-2ts}\,|\widehat h(\xi)|^2\dd\xi}
       {\int_{\R^d}|\widehat h(\xi)|^2\dd\xi}.
\]
The supremum of this quotient over all Schwartz functions $h\ne0$ is the
supremum of the continuous function $s^{2m}e^{-2ts}$ over $s\ge0$, as one
sees by choosing $\widehat h$ concentrated near a point $\xi$ at which the
supremum is approached. For $m\ge1$ this function attains its maximum at
$s=m/t$, where it equals $(m/t)^{2m}e^{-2m}$. Since the left-hand side is
at most $\opnm{\A^mP_t}^2$, we conclude that
\[
 \opnm{\A^mP_t}^2\ge\sup_{s\ge0}s^{2m}e^{-2ts}
       =\left(\frac m{et}\right)^{2m},
\]
which is \eqref{eq:sharp-even}.

For the gradient we use that, for a function of the form
$f(x,v)=g(x-v/\gamma)$, one has $\nabla_xf=(\nabla g)(x-v/\gamma)$ and
$\nabla_vf=-\gamma^{-1}(\nabla g)(x-v/\gamma)$, so that
\[
 |\nabla_{x,v}f|^2=(1+\gamma^{-2})\,|\nabla g|^2(x-v/\gamma)
 \qquad\text{and}\qquad
 \nm{\nabla_{x,v}\,\mathcal Ig}^2
 =(1+\gamma^{-2})\,\|\nabla g\|_{L^2(\dd x)}^2 ,
\]
by the same computation as for the isometry property. We apply this to the
function $g$ of \eqref{eq:heatpowers}. Its gradient has the Fourier
transform $i\xi\,(-s)^me^{-ts}\,\widehat h(\xi)$, and
$|\xi|^2=(\gamma/\sigma)s$, so that
\[
 \frac{\nm{\nabla_{x,v}\A^mP_t\,\mathcal Ih}^2}{\nm{\mathcal Ih}^2}
 =(1+\gamma^{-2})\,\frac\gamma\sigma\,
  \frac{\int_{\R^d}s^{2m+1}e^{-2ts}\,|\widehat h(\xi)|^2\dd\xi}
       {\int_{\R^d}|\widehat h(\xi)|^2\dd\xi}.
\]
As before, the supremum over $h$ is the supremum of $s^{2m+1}e^{-2ts}$
over $s\ge0$, which is attained at $s=(2m+1)/(2t)$ and equals
$((2m+1)/(2et))^{2m+1}$. Hence, for every $m\ge0$,
\[
 \sup_{h\ne0}\frac{\nm{\nabla_{x,v}\A^mP_th}^2}{\nm h^2}
 \ge\frac{\gamma}{\sigma}(1+\gamma^{-2})
       \left(\frac{2m+1}{2et}\right)^{2m+1},
\]
which is \eqref{eq:sharp-odd}. Finally, if the power $t^{-2m}$ in
\eqref{eq:main-even} or the power $t^{-(2m+1)}$ in \eqref{eq:main-odd}
could be replaced by a larger power of $1/t$, uniformly over a class of
potentials containing $L=0$, this would contradict the lower bounds just
proved for large $t$.
\end{proof}

\section*{Statement on the use of AI}

During the preparation of this work we used Claude to assist
with the search of the literature, with the wording
of the manuscript, with the drawing of the figures and with the
derivation of the explicit values of the constants. The example of
Remark~\ref{rem:CL} was also obtained with
its assistance. After using this tool, we reviewed and edited the content as
needed, and we take full responsibility for the content of the article.

\section*{Data availability}

Data sharing is not applicable to this article, as no datasets were
generated or analysed during the current study.

\section*{Conflict of interest}

The author has no competing interests to declare that are relevant to the
content of this article.

\appendix

\section{Proofs of Propositions~\ref{prop:maxdissipative}
and~\ref{prop:adjoint}}
\label{apx:maxdissipative}

This appendix contains the proofs of the two propositions of
Section~\ref{sec:generator}, both of which are known. The first one is
obtained by transporting to $L^2(\mu)$, through a unitary transformation, a
result which is stated in the literature for the unweighted space
$L^2(\dd x\dd v)$, and the second one is a short argument based on the
first.

\subsection{Proof of Proposition~\ref{prop:maxdissipative}}

The literature on kinetic Fokker--Planck operators uses the sign
convention opposite to ours: its statements concern $-\A$ rather than
$\A$. Translated into our convention, the result of Helffer and Nier
\cite[Proposition~5.5]{HelfferNier} is the following: for every smooth
potential $V$, the differential operator
\[
 -\Delta_v+\frac{|v|^2}{4}-\frac d2+v\cdot\nabla_x-\nabla V\cdot\nabla_v,
\]
initially defined on $\Cc$, is closable in $L^2(\dd x\dd v)$, and the
negative of its closure is maximal dissipative; the same statement with a
friction parameter is
\cite[Proposition~3.3]{NierLectures}. We reduce $-\A$ to a positive multiple
of this operator by two unitary transformations.

The first one is the multiplication by $e^{-(L+|v|^2/2)/(2\sigma)}$, which
maps $L^2(\mu)$ isometrically onto $L^2(\dd x\dd v)$, also when
$\mu(\R^{2d})=\infty$. Since $\Lham(L+|v|^2/2)=0$, the operator $\Lham$
commutes with this multiplication, and a direct computation of the
conjugate of $\LFD$ gives, on $\Cc$,
\begin{equation}\label{eq:KFPconj}
 -e^{-(L+|v|^2/2)/(2\sigma)}\,\A\, e^{(L+|v|^2/2)/(2\sigma)}
 =\Lham-\gamma\sigma\Delta_v
                  +\frac{\gamma}{4\sigma}|v|^2-\frac{\gamma d}{2}.
\end{equation}
The second one is the dilation $x=(\sqrt\sigma/\gamma)\,x'$,
$v=\sqrt\sigma\,v'$, normalized so as to be unitary, which transforms the
right-hand side of \eqref{eq:KFPconj} into
\begin{equation}\label{eq:HNoperator}
 \gamma\left(-\Delta_{v'}+\frac{|v'|^2}{4}-\frac d2
 +v'\cdot\nabla_{x'}-\nabla V(x')\cdot\nabla_{v'}\right),
 \qquad V(x')=\frac1\sigma\,L\Big(\frac{\sqrt\sigma}{\gamma}\,x'\Big),
\end{equation}
that is, into $\gamma$ times the operator of the cited result for the
smooth potential $V$. Let $U$ denote the composition of the two
transformations. It is a unitary operator from $L^2(\mu)$ onto
$L^2(\dd x'\dd v')$ which maps $\Cc$ onto the test functions of the new
variables, and $U(-\A)U^{-1}$ coincides with \eqref{eq:HNoperator} on test
functions. Consequently the closures of the two operators are unitarily
equivalent through $U$. Since unitary equivalence and multiplication by a
positive constant preserve maximal dissipativity, the closure of $\A$ is
maximal dissipative.

We turn to the description \eqref{eq:maxdomain} of the domain. On its
right-hand side, $\A h$ is computed in the sense of distributions and is
required to belong to $L^2(\mu)$; equivalently, there is $f\in L^2(\mu)$
such that $\ip h{(\Lham+\gamma\LFD)\varphi}=\ip f\varphi$ for every
$\varphi\in\Cc$. This is the largest subspace of $L^2(\mu)$ on which the
differential expression $\A$ defines an operator with values in $L^2(\mu)$,
whence its name, and \eqref{eq:maxdomain} says that it coincides with the
domain of the closure. For the operator \eqref{eq:HNoperator} this is
stated in \cite[Remark~3.2]{NierLectures}, and it is transported to $\A$ by
$U^{-1}$, which acts on distributions as well. Since that remark is stated
without proof, an independent derivation of \eqref{eq:maxdomain} is
included in the proof of Proposition~\ref{prop:adjoint}.

Finally, a maximal dissipative operator generates a strongly continuous
semigroup of contractions by the Lumer--Phillips theorem
\cite[Theorem~II.3.15]{EngelNagel}, and the fact that $h_t=P_th_0$ is the
unique solution of \eqref{eq:A} for $h_0\in\Dom(\A)$ is
\cite[Proposition~II.6.2]{EngelNagel}. This completes the proof of
Proposition~\ref{prop:maxdissipative}.

\subsection{Proof of Proposition~\ref{prop:adjoint}}

Denote by $\A'$ the closure of $\Lham+\gamma\LFD$ on $\Cc$. Velocity
reversal preserves $\Cc$ and, by \eqref{eq:Jcore}, exchanges the operators
$\A$ and $\Lham+\gamma\LFD$ on $\Cc$; since $\J$ is unitary, taking closures gives $\A'=\J\A\J$ with
$\Dom(\A')=\J\Dom(\A)$, and $\A'$ is maximal dissipative because $\A$ is,
by Proposition~\ref{prop:maxdissipative}.

We next compare $\A'$ with $\A^*$. Let $h\in\Dom(\A)$, and let
$h_n\in\Cc$ be such that $h_n\to h$ and $\A h_n\to\A h$ in $L^2(\mu)$. For
$\varphi\in\Cc$, the formal adjointness \eqref{eq:formaladjoint} gives
\[
 \ip{\A h}{\varphi}=\lim_{n\to\infty}\ip{\A h_n}{\varphi}
   =\lim_{n\to\infty}\ip{h_n}{(\Lham+\gamma\LFD)\varphi}
   =\ip h{(\Lham+\gamma\LFD)\varphi}.
\]
Hence $\Cc\subset\Dom(\A^*)$ and $\A^*$ extends $\Lham+\gamma\LFD$ on
$\Cc$; since adjoints are closed, $\A^*$ extends $\A'$. On the other hand
$\A^*$ is dissipative: for $\lambda>0$ the operator
$\lambda-\A^*=(\lambda-\A)^*$ is invertible with
$\opnm{(\lambda-\A^*)^{-1}}=\opnm{(\lambda-\A)^{-1}}\le1/\lambda$, by
\eqref{eq:halfplane}, so that $\nm{(\lambda-\A^*)h}\ge\lambda\nm h$ for
$h\in\Dom(\A^*)$; squaring this inequality and letting $\lambda\to\infty$
gives $\Ree\ip{\A^*h}{h}\le0$. Thus $\A^*$ is a
dissipative extension of the maximal dissipative operator $\A'$, and the
two operators coincide. This proves \eqref{eq:J}.

The same argument yields the maximal-domain characterization
\eqref{eq:maxdomain}. Indeed, $\A$ is closed and densely defined, so that
$\A=\A^{**}=(\A')^*$, and the adjoint of a closure coincides with the
adjoint of the operator before closure. Hence $\A$ is the adjoint of the
operator $\Lham+\gamma\LFD$ with domain $\Cc$, and by the very definition
of the adjoint its domain consists of those $h\in L^2(\mu)$ for which there
is $f\in L^2(\mu)$ with $\ip h{(\Lham+\gamma\LFD)\varphi}=\ip f\varphi$ for
all $\varphi\in\Cc$, that is, for which $\A h$, computed in the sense of
distributions, belongs to $L^2(\mu)$. This is precisely the right-hand side
of \eqref{eq:maxdomain}.

We turn to the consequences. Since $\J$ is bounded and selfadjoint,
$(\J\A)^*=\A^*\J=\J\A\J\J=\J\A$, and the domain of $\A^*\J$ is
$\{h:\J h\in\Dom(\A^*)\}=\Dom(\A)$ by \eqref{eq:J}; thus $\J\A$ is
selfadjoint on $\Dom(\A)$. In particular
$\ip{\A h}{\J h}=\ip{\J\A h}{h}$ is real, which is \eqref{eq:Jreal}.
Finally, let $h\in\Dom(\A)$ and let $h_n\in\Cc$ be such that $h_n\to h$
and $\A h_n\to\A h$ in $L^2(\mu)$. By \eqref{eq:dirichlet}
applied to differences,
\[
 \gamma\sigma\nm{\nabla_v(h_n-h_j)}^2
       =-\Ree\ip{\A(h_n-h_j)}{h_n-h_j}\longrightarrow0
 \qquad(n,j\to\infty),
\]
so that $\nabla_vh_n$ converges in $L^2(\mu)$, necessarily to the weak
derivative $\nabla_vh$, and \eqref{eq:dirichlet} passes to the limit.
This completes the proof of Proposition~\ref{prop:adjoint}.

\section{Proofs of Lemmas~\ref{lem:reverse} and~\ref{lem:static}}
\label{apx:energy}

This appendix contains the proofs of the two energy inequalities of
Section~\ref{sec:energy}. Throughout this appendix we assume that
Assumption~\ref{ass:WH} holds.
Since test functions are dense in
$H^1_x(\mu)$, the inequality \eqref{eq:WH} extends to all of this space by
Fatou's lemma, and applying it to the components of $\nabla_vh$ and of
$\nabla_xh$ gives the two relative bounds
\begin{equation}\label{eq:relative}
 \begin{aligned}
 \nm{(\nabla^2L)\nabla_vh}^2
 &\le C_L\big(\nm{\nabla_vh}^2+\nm{\nabla_x\nabla_vh}^2\big),\\
 \nm{(\nabla^2L)\nabla_xh}^2
 &\le C_L\big(\nm{\nabla_xh}^2+\nm{\nabla_x^2h}^2\big),
 \end{aligned}
\end{equation}
whenever the right-hand sides are finite; these are the only forms in which
the assumption is used in this appendix and in the next one.

The two energy inequalities, Lemmas~\ref{lem:reverse}
and~\ref{lem:static}, are known when $\mu$ is a probability measure. The first one rests on the time-weighted
functional \eqref{eq:firstenergy}, which was introduced by H\'erau
\cite[Theorem~1.1 and Section~3]{Herau} for potentials whose derivatives of
order at least two are bounded. Zhang
\cite[Theorems~2 and~3, Proposition~16]{Zhang} proved that this functional
is nonincreasing along the flow under an assumption equivalent to
Assumption~\ref{ass:WH}, in the case where $\gamma=\sigma=1$ and $\mu$ is a
probability measure; his proof ends with the inequality~(5.5), which is
equivalent to \eqref{eq:reverse-t}. The functional $\cG$ of the second
lemma is the modified $H^1$ norm of Villani \cite[Theorem~18]{Villani},
and \eqref{eq:static} is his inequality~(4.1). The proofs given
here are not different in essence from those of the
finite-mass case: the differentiation of the functionals along the flow is
the computation of the references just quoted, with the parameters
$\gamma$ and $\sigma$, and it uses nowhere that $\mu$ is finite. What we
add is the justification of this computation for initial data in
$L^2(\mu)$, by means of an elliptic regularization; this step is called a
standard approximation in \cite{Zhang} and is not written out there, nor
in \cite[Appendix~A.21.2]{Villani}, where the case of infinite mass is
described under the pointwise condition \eqref{eq:VG}.

\subsection{The elliptic regularization}\label{sec:ellipticreg}

The energy identities below involve derivatives of $P_th$ of order up to
two, and for $L^2(\mu)$ data their existence is part of what has to be
proved. The corresponding approximation argument is not carried out in the
references (like \cite[Appendix~A.21]{Villani} and \cite{Zhang}) on which
the energy estimates of Section~\ref{sec:energy} are based, and we
therefore give it in full. To justify the computations we regularize the generator by adding a
small multiple of the symmetric position operator
\[
 \Sx=\sigma\Delta_x-\nabla L\cdot\nabla_x ,
 \qquad \ip{\Sx f}g=-\sigma\ip{\nabla_xf}{\nabla_xg}
 \quad(f,g\in\Cc),
\]
which preserves the weight $\mu$. For $0<\varepsilon\le1$ we set
\begin{equation}\label{eq:ellipticA}
 \A_\varepsilon=\A+\varepsilon\Sx,
 \qquad P_t^\varepsilon=e^{t\A_\varepsilon}.
\end{equation}
The properties of this regularization which we shall use are collected in
the next lemma, whose proof is given in Appendix~\ref{apx:elliptic}.

\begin{lemma}[Elliptic regularization]\label{lem:elliptic}
Let $0<\varepsilon\le1$.
\begin{enumerate}[label=\textup{(\roman*)},leftmargin=*]
\item The closure of $\A_\varepsilon$ on $\Cc$ is maximal dissipative with
maximal domain, and it generates an analytic contraction semigroup
$(P_t^\varepsilon)_{t\ge0}$ on $L^2(\mu)$. For
$f,g\in\Dom(\A_\varepsilon)$,
\begin{equation}\label{eq:ellipticdiss}
 \Ree\big(\ip{\A_\varepsilon f}{g}+\ip f{\A_\varepsilon g}\big)
 =-2\gamma\sigma\Ree\ip{\nabla_vf}{\nabla_vg}
  -2\varepsilon\sigma\Ree\ip{\nabla_xf}{\nabla_xg}.
\end{equation}
\item The domain of $\A_\varepsilon$ is contained in $H^2(\mu)$, and there
is a constant $C_\varepsilon$ such that
\begin{equation}\label{eq:ellipticgraph}
 \|h\|_{H^2(\mu)}
       \le C_\varepsilon(\nm{\A_\varepsilon h}+\nm h),
       \qquad h\in\Dom(\A_\varepsilon).
\end{equation}
\item For every $h\in L^2(\mu)$ and every $t_0>0$,
\begin{equation}\label{eq:visclimit}
 \sup_{0\le t\le t_0}\nm{P_t^\varepsilon h-P_th}\longrightarrow0
 \qquad(\varepsilon\to0).
\end{equation}
\item For $h\in\Cc$ the trajectory $u_t=P_t^\varepsilon h$ belongs to
$C^1([0,\infty);H^3(\mu))$, and for all multi-indices $\alpha,\beta$ in
the variables $(x,v)$ of order at most one the functions
$\partial^\alpha u_t$ belong to $\Dom(\A_\varepsilon)$, the commutators
$[\partial^\alpha,\A]u_t$ and $[\partial^\alpha,\Sx]u_t$ belong to
$L^2(\mu)$, and
\begin{equation}\label{eq:star}
 \begin{aligned}
 \frac{\dd}{\dd t}\Ree\ip{\partial^\alpha u_t}{\partial^\beta u_t}
 ={}&-2\gamma\sigma\Ree
      \ip{\nabla_v\partial^\alpha u_t}{\nabla_v\partial^\beta u_t}
    -2\varepsilon\sigma\Ree
      \ip{\nabla_x\partial^\alpha u_t}{\nabla_x\partial^\beta u_t}\\
 &+\Ree\ip{[\partial^\alpha,\A_\varepsilon]u_t}{\partial^\beta u_t}
  +\Ree\ip{\partial^\alpha u_t}{[\partial^\beta,\A_\varepsilon]u_t}.
 \end{aligned}
\end{equation}
\end{enumerate}
\end{lemma}

Since $[\partial^\alpha,\A_\varepsilon]=[\partial^\alpha,\A]
+\varepsilon[\partial^\alpha,\Sx]$, the right-hand side of \eqref{eq:star}
splits into the terms which do not carry an explicit factor $\varepsilon$
and those which do. Consider a quadratic expression of the form
\begin{equation}\label{eq:quadexpr}
 \Psi(t;w)=\sum_ic_i\,t^{p_i}\,\Ree\ip{\partial^{\alpha_i}w}{\partial^{\beta_i}w},
\end{equation}
where the sum is finite, the $c_i$ are real, the $p_i$ are nonnegative
integers and the multi-indices have order at most one. For a regularized
trajectory $u_t=P_t^\varepsilon h$ with $h\in\Cc$ we define its
\emph{kinetic part} and its \emph{viscous part} by
\begin{align}
 \Kin[\Psi](t)&=\sum_ic_ip_i\,t^{p_i-1}
     \Ree\ip{\partial^{\alpha_i}u_t}{\partial^{\beta_i}u_t}\notag\\
  &\quad+\sum_ic_it^{p_i}\Big\{-2\gamma\sigma\Ree
      \ip{\nabla_v\partial^{\alpha_i}u_t}{\nabla_v\partial^{\beta_i}u_t}\notag\\
  &\qquad\qquad\qquad
   +\Ree\ip{[\partial^{\alpha_i},\A]u_t}{\partial^{\beta_i}u_t}
   +\Ree\ip{\partial^{\alpha_i}u_t}{[\partial^{\beta_i},\A]u_t}\Big\},
   \label{eq:kinpart}\\
 \Vis[\Psi](t)&=\sum_ic_it^{p_i}\Big\{-2\sigma\Ree
      \ip{\nabla_x\partial^{\alpha_i}u_t}{\nabla_x\partial^{\beta_i}u_t}\notag\\
  &\qquad\qquad\qquad
   +\Ree\ip{[\partial^{\alpha_i},\Sx]u_t}{\partial^{\beta_i}u_t}
   +\Ree\ip{\partial^{\alpha_i}u_t}{[\partial^{\beta_i},\Sx]u_t}\Big\},
   \label{eq:vispart}
\end{align}
so that, by \eqref{eq:star},
\begin{equation}\label{eq:split}
 \frac{\dd}{\dd t}\Psi(t;u_t)=\Kin[\Psi](t)+\varepsilon\,\Vis[\Psi](t).
\end{equation}
The kinetic part is the expression obtained by differentiating
$\Psi(t;P_th)$ formally along the kinetic semigroup itself.

\subsection{Gradient identities and the choice of coefficients}

The commutators of the first derivatives with $\A$ and with $\Sx$ are
\[
 [\nabla_v,\A]=-\nabla_x-\gamma\nabla_v,\qquad
 [\nabla_x,\A]=(\nabla^2L)\nabla_v,\qquad
 [\nabla_v,\Sx]=0,\qquad
 [\nabla_x,\Sx]=-(\nabla^2L)\nabla_x .
\]
Inserting them into \eqref{eq:kinpart}, we obtain the kinetic parts of the
four basic quadratic quantities,
\begin{align}
 \Kin\big[\nm{u_t}^2\big]
 &=-2\gamma\sigma\nm{\nabla_vu_t}^2,\label{eq:baseE}\\
 \Kin\big[\nm{\nabla_vu_t}^2\big]
 &=-2\gamma\sigma\nm{\nabla_v^2u_t}^2
   -2\gamma\nm{\nabla_vu_t}^2
   -2\Ree\ip{\nabla_vu_t}{\nabla_xu_t},\label{eq:EU}\\
 \Kin\big[\nm{\nabla_xu_t}^2\big]
 &=-2\gamma\sigma\nm{\nabla_v\nabla_xu_t}^2
   +2\Ree\ip{\nabla_xu_t}{(\nabla^2L)\nabla_vu_t},\label{eq:EW}\\
 \Kin\big[\Ree\ip{\nabla_vu_t}{\nabla_xu_t}\big]
 &=-2\gamma\sigma\Ree\ip{\nabla_v^2u_t}{\nabla_v\nabla_xu_t}
   -\nm{\nabla_xu_t}^2\notag\\
 &\quad-\gamma\Ree\ip{\nabla_vu_t}{\nabla_xu_t}
   +\Ree\ip{\nabla_vu_t}{(\nabla^2L)\nabla_vu_t}.
                                                  \label{eq:EC}
\end{align}
The mixed quantity in \eqref{eq:EC} is the one which produces the
dissipation $-\nm{\nabla_xu_t}^2$ of the position gradient. Similarly, \eqref{eq:vispart}
gives the viscous parts
\begin{equation}\label{eq:visparts}
 \begin{aligned}
 \Vis\big[\nm{u_t}^2\big]&=-2\sigma\nm{\nabla_xu_t}^2,\qquad
 \Vis\big[\nm{\nabla_vu_t}^2\big]=-2\sigma\nm{\nabla_x\nabla_vu_t}^2,\\
 \Vis\big[\nm{\nabla_xu_t}^2\big]&=-2\sigma\nm{\nabla_x^2u_t}^2
     -2\Ree\ip{(\nabla^2L)\nabla_xu_t}{\nabla_xu_t},\\
 \Vis\big[\Ree\ip{\nabla_vu_t}{\nabla_xu_t}\big]
   &=-2\sigma\Ree\ip{\nabla_x\nabla_vu_t}{\nabla_x^2u_t}
     -\Ree\ip{(\nabla^2L)\nabla_vu_t}{\nabla_xu_t}.
 \end{aligned}
\end{equation}

We now fix, once and for all, three positive coefficients $a,b,c$
satisfying
\begin{equation}\label{eq:coeff}
 \begin{gathered}
 ac=4b^2,\qquad c\le\frac b6,\qquad c\le\frac{\gamma\sigma b}{16C_L},\\
 a+2b\sqrt{C_L}+\frac{4b^2C_L}{\gamma\sigma c}
       +\frac{4(a+2b+b\gamma+c\sqrt{C_L})^2}{b}\le\gamma\sigma.
 \end{gathered}
\end{equation}
An admissible choice is
\begin{equation}\label{eq:coeffexplicit}
\begin{aligned}
 \frac cb&=\min\left\{\frac16,\frac{\gamma\sigma}{16(1+C_L)}\right\},\qquad
 \frac ab=\frac4{c/b},\\
 b&=\min\left\{1,\frac{\gamma\sigma}{1+a/b+2\sqrt{C_L}
   +\dfrac{4C_L}{\gamma\sigma(c/b)}
   +4\big(a/b+2+\gamma+(c/b)\sqrt{C_L}\big)^2}\right\}.
\end{aligned}
\end{equation}
The relation $ac=4b^2$ makes the quadratic form with coefficients
$a,2b,c$ comparable to its diagonal part: for all complex numbers, or
vectors, $X$ and $Y$,
\begin{equation}\label{eq:form-positive}
 \tfrac12\big(a|X|^2+c|Y|^2\big)
 \le a|X|^2+2b\Ree(X\overline Y)+c|Y|^2
 \le\tfrac32\big(a|X|^2+c|Y|^2\big),
\end{equation}
because $2b|X||Y|=\sqrt{ac}\,|X||Y|\le\tfrac12(a|X|^2+c|Y|^2)$. The same
holds with $a,b,c$ replaced by $at,bt^2,ct^3$ for any $t>0$, since
$(at)(ct^3)=4(bt^2)^2$.

With these coefficients we introduce, for $w\in H^1(\mu)$ and $t\ge0$, the
time-weighted energy
\begin{equation}\label{eq:firstenergy}
 \cF_t(w)=\nm{w}^2+at\nm{\nabla_vw}^2
 +2bt^2\Ree\ip{\nabla_vw}{\nabla_xw}
 +ct^3\nm{\nabla_xw}^2,
\end{equation}
so that $\cF_1$ is the energy $\cG$ defined in \eqref{eq:staticenergy}. By
\eqref{eq:form-positive}, integrated against $\mu$,
\begin{equation}\label{eq:Fequiv}
 \nm w^2+\tfrac12\big(at\nm{\nabla_vw}^2+ct^3\nm{\nabla_xw}^2\big)
 \le\cF_t(w)\le
 \nm w^2+\tfrac32\big(at\nm{\nabla_vw}^2+ct^3\nm{\nabla_xw}^2\big),
\end{equation}
and the case $t=1$ is \eqref{eq:Gequiv}.

\subsection{Proof of Lemma~\ref{lem:reverse}}

\emph{Step 1: the regularized trajectories.} Let $h\in\Cc$, let
$0<\varepsilon\le1$ and let $u_t=P_t^\varepsilon h$. By
Lemma~\ref{lem:elliptic}(iv), the function
$t\mapsto\cF_t(u_t)$ is continuously differentiable on $[0,1]$, it equals
$\nm h^2$ at $t=0$, and by \eqref{eq:split}
\[
 \frac{\dd}{\dd t}\cF_t(u_t)=\Kin[\cF](t)+\varepsilon\Vis[\cF](t).
\]

\emph{Step 2: the kinetic part.} By \eqref{eq:baseE}--\eqref{eq:EC}, the
kinetic part is
\begin{align*}
 \Kin[\cF](t)={}&(-2\gamma\sigma+a)\nm{\nabla_vu_t}^2
 +(3c-2b)t^2\nm{\nabla_xu_t}^2\\
 &+(4b-2a)t\Ree\ip{\nabla_vu_t}{\nabla_xu_t}
 -2b\gamma t^2\Ree\ip{\nabla_vu_t}{\nabla_xu_t}
 -2a\gamma t\nm{\nabla_vu_t}^2\\
 &+2bt^2\Ree\ip{\nabla_vu_t}{(\nabla^2L)\nabla_vu_t}
 +2ct^3\Ree\ip{\nabla_xu_t}{(\nabla^2L)\nabla_vu_t}\\
 &-2\gamma\sigma\big(
 at\nm{\nabla_v^2u_t}^2+
 2bt^2\Ree\ip{\nabla_v^2u_t}{\nabla_v\nabla_xu_t}
 +ct^3\nm{\nabla_v\nabla_xu_t}^2\big).
\end{align*}
By \eqref{eq:form-positive} the last bracket is at least
$\tfrac12(at\nm{\nabla_v^2u_t}^2+ct^3\nm{\nabla_v\nabla_xu_t}^2)$. In the two
terms containing the Hessian we use the relative bound
\eqref{eq:relative} in the form
$\nm{(\nabla^2L)\nabla_vu_t}\le\sqrt{C_L}\,(\nm{\nabla_vu_t}
+\nm{\nabla_v\nabla_xu_t})$. Since $t\le1$, the terms which may be positive are then controlled by the following three
applications of Young's inequality:
\begin{align*}
 2b\sqrt{C_L}\,t^2\nm{\nabla_vu_t}\nm{\nabla_v\nabla_xu_t}
 &\le\frac{\gamma\sigma c}{4}t^3\nm{\nabla_v\nabla_xu_t}^2
       +\frac{4b^2C_L}{\gamma\sigma c}\nm{\nabla_vu_t}^2,\\
 2c\sqrt{C_L}\,t^3\nm{\nabla_xu_t}\nm{\nabla_v\nabla_xu_t}
 &\le\frac b4t^2\nm{\nabla_xu_t}^2
       +\frac{\gamma\sigma c}{4}t^3\nm{\nabla_v\nabla_xu_t}^2,\\
 2(a+2b+b\gamma+c\sqrt{C_L})\,t\,\nm{\nabla_vu_t}\nm{\nabla_xu_t}
 &\le\frac b4t^2\nm{\nabla_xu_t}^2
  +\frac{4(a+2b+b\gamma+c\sqrt{C_L})^2}{b}
                                    \nm{\nabla_vu_t}^2.
\end{align*}
The second inequality uses $c\le\gamma\sigma b/(16C_L)$, and the third one
collects all the cross terms between $\nabla_vu_t$ and $\nabla_xu_t$, namely
those with coefficients $4b-2a$ and $2b\gamma$ and the lower-order part of
the last Hessian term. The remaining diagonal Hessian contribution is at
most $2b\sqrt{C_L}\nm{\nabla_vu_t}^2$. Collecting the coefficients, and using
$c\le b/6$ for the coefficient of $t^2\nm{\nabla_xu_t}^2$ and the last
inequality in \eqref{eq:coeff} for the coefficient of $\nm{\nabla_vu_t}^2$,
we conclude that
\begin{equation}\label{eq:F-diss}
 \Kin[\cF](t)\le
 -\gamma\sigma\nm{\nabla_vu_t}^2
 -bt^2\nm{\nabla_xu_t}^2
 -\gamma\sigma at\nm{\nabla_v^2u_t}^2
 -\frac{\gamma\sigma c}{2}t^3\nm{\nabla_v\nabla_xu_t}^2\le0,
 \qquad0<t\le1.
\end{equation}

\emph{Step 3: the viscous part.} By \eqref{eq:visparts},
\begin{equation}\label{eq:viscerror}
 \begin{aligned}
 \Vis[\cF](t)={}&-2\sigma\big(
 \nm{\nabla_xu_t}^2+
 at\nm{\nabla_x\nabla_vu_t}^2
 +2bt^2\Ree
       \ip{\nabla_x\nabla_vu_t}{\nabla_x^2u_t}
                  +ct^3\nm{\nabla_x^2u_t}^2\big)\\
 &-2bt^2\Ree
       \ip{(\nabla^2L)\nabla_vu_t}{\nabla_xu_t}
       -2ct^3\Ree
       \ip{(\nabla^2L)\nabla_xu_t}{\nabla_xu_t}.
 \end{aligned}
\end{equation}
The first bracket is nonnegative, and by \eqref{eq:form-positive} it is at
least $\nm{\nabla_xu_t}^2+\tfrac12(at\nm{\nabla_x\nabla_vu_t}^2
+ct^3\nm{\nabla_x^2u_t}^2)$. The two force terms are estimated with both
inequalities in \eqref{eq:relative}:
\begin{align*}
 &2b\sqrt{C_L}\,t^2
                    \nm{\nabla_vu_t}\nm{\nabla_xu_t}
 \le C\big(at\nm{\nabla_vu_t}^2+
                           ct^3\nm{\nabla_xu_t}^2\big),\\
 &2b\sqrt{C_L}\,t^2
            \nm{\nabla_x\nabla_vu_t}\nm{\nabla_xu_t}
 \le\frac{\sigma a}{4}t
                   \nm{\nabla_x\nabla_vu_t}^2
       +\frac{4b^2C_L}{\sigma a}t^3
                                     \nm{\nabla_xu_t}^2,\\
 &2c\sqrt{C_L}\,t^3
                      \nm{\nabla_x^2u_t}\nm{\nabla_xu_t}
 \le\frac{\sigma c}{4}t^3
                         \nm{\nabla_x^2u_t}^2
       +\frac{4cC_L}{\sigma}t^3\nm{\nabla_xu_t}^2,
\end{align*}
together with the term $2c\sqrt{C_L}\,t^3\nm{\nabla_xu_t}^2$. The two
second-derivative terms are absorbed by the first bracket of
\eqref{eq:viscerror}, and the remaining terms are bounded by a constant
times $at\nm{\nabla_vu_t}^2+ct^3\nm{\nabla_xu_t}^2$, hence by a constant times
$\cF_t(u_t)$ in view of \eqref{eq:Fequiv}. Thus
$\Vis[\cF](t)\le C\cF_t(u_t)$ with a constant $C$ depending only on
$a,b,c,\sigma$ and $C_L$, and in particular independent of $\varepsilon$.

\emph{Step 4: Gronwall's inequality.} Steps~2 and~3 give
\begin{equation}\label{eq:viscF}
 \frac{\dd}{\dd t}\cF_t(u_t)\le C\varepsilon\,\cF_t(u_t),\qquad0<t\le1,
 \qquad\text{and hence}\qquad
 \cF_t(P_t^\varepsilon h)\le e^{C\varepsilon t}\nm h^2 .
\end{equation}

\emph{Step 5: the limit $\varepsilon\to0$.} Fix $t\in(0,1]$. By
\eqref{eq:visclimit}, $P_t^\varepsilon h\to P_th$ in $L^2(\mu)$, and by
\eqref{eq:viscF} and \eqref{eq:Fequiv} the gradients
$\nabla_{x,v}P_t^\varepsilon h$ are bounded in $L^2(\mu)$ uniformly in
$\varepsilon$. Along a sequence $\varepsilon_n\to0$ they converge weakly in
$L^2(\mu)$, and the limit is the distributional gradient of $P_th$, because
distributional derivatives pass to the limit when the weight is locally
smooth and positive. Hence $P_th\in H^1(\mu)$. The functional $\cF_t$ is
the sum of $\nm{\cdot}^2$ and of a nonnegative continuous quadratic form in
the gradient, and is therefore lower semicontinuous for this mode of
convergence; letting $n\to\infty$ in \eqref{eq:viscF} we obtain
\[
 \cF_t(P_th)\le\nm h^2,\qquad h\in\Cc,\quad0<t\le1.
\]

\emph{Step 6: general initial data.} For $h\in L^2(\mu)$ we choose
$h_n\in\Cc$ converging to $h$ in $L^2(\mu)$. Then $P_th_n\to P_th$ in
$L^2(\mu)$ and $\cF_t(P_th_n)\le\nm{h_n}^2$, and the same weak compactness
and lower semicontinuity argument gives $P_th\in H^1(\mu)$ and
$\cF_t(P_th)\le\nm h^2$. The lower bound in \eqref{eq:Fequiv} now yields
\eqref{eq:reverse-t}, and \eqref{eq:reverse-one} is the case $t=1$.
This completes the proof of Lemma~\ref{lem:reverse}. \qed

\subsection{Proof of Lemma~\ref{lem:static}}

Let first $h\in\Cc$ and $u_t=P_t^\varepsilon h$ with $0<\varepsilon\le1$.
By \eqref{eq:split},
$\frac{\dd}{\dd t}\cG(u_t)=\Kin[\cG](t)+\varepsilon\Vis[\cG](t)$, where
this time the coefficients of the quadratic expression do not depend on
$t$. The gradient identities \eqref{eq:baseE}--\eqref{eq:EC} give
\begin{align*}
 \Kin[\cG](t)={}&-2\gamma\sigma\nm{\nabla_vu_t}^2
   -2a\gamma\nm{\nabla_vu_t}^2-2b\nm{\nabla_xu_t}^2
   -2(a+b\gamma)\Ree\ip{\nabla_vu_t}{\nabla_xu_t}\\
 &+2b\Ree\ip{\nabla_vu_t}{(\nabla^2L)\nabla_vu_t}
   +2c\Ree\ip{\nabla_xu_t}{(\nabla^2L)\nabla_vu_t}\\
 &-2\gamma\sigma\big(a\nm{\nabla_v^2u_t}^2
   +2b\Ree\ip{\nabla_v^2u_t}{\nabla_v\nabla_xu_t}
   +c\nm{\nabla_v\nabla_xu_t}^2\big).
\end{align*}
The last bracket is bounded below by \eqref{eq:form-positive}, the Hessian
terms are bounded by \eqref{eq:relative}, and the three required Young
inequalities are
\begin{align*}
2b\sqrt{C_L}\nm{\nabla_vu_t}\nm{\nabla_v\nabla_xu_t}
&\le\frac{\gamma\sigma c}{4}\nm{\nabla_v\nabla_xu_t}^2
 +\frac{4b^2C_L}{\gamma\sigma c}\nm{\nabla_vu_t}^2,\\
2c\sqrt{C_L}\nm{\nabla_xu_t}\nm{\nabla_v\nabla_xu_t}
&\le\frac b2\nm{\nabla_xu_t}^2+
                 \frac{2c^2C_L}{b}\nm{\nabla_v\nabla_xu_t}^2,\\
2(a+b\gamma+c\sqrt{C_L})\nm{\nabla_vu_t}\nm{\nabla_xu_t}
&\le\frac b2\nm{\nabla_xu_t}^2+
       \frac{2(a+b\gamma+c\sqrt{C_L})^2}{b}\nm{\nabla_vu_t}^2.
\end{align*}
Since $2c^2C_L/b\le\gamma\sigma c/8$ by \eqref{eq:coeff}, and since the
coefficient of $\nm{\nabla_vu_t}^2$ is controlled by the last inequality in
\eqref{eq:coeff} as before, we obtain the analogue of \eqref{eq:F-diss}
with all powers of $t$ replaced by one,
\[
 \Kin[\cG](t)\le-\gamma\sigma\nm{\nabla_vu_t}^2-b\nm{\nabla_xu_t}^2
   -\gamma\sigma a\nm{\nabla_v^2u_t}^2
   -\frac{\gamma\sigma c}{2}\nm{\nabla_v\nabla_xu_t}^2
   \le-c_*\,\cD(u_t).
\]
The viscous part is given by \eqref{eq:viscerror} with all powers of $t$
replaced by one, and the estimates of Step~3 of the proof of
Lemma~\ref{lem:reverse}, with constant time weights, give
$\Vis[\cG](t)\le C\cG(u_t)$. Therefore
\[
 \frac{\dd}{\dd t}\cG(u_t)\le-c_*\cD(u_t)+C\varepsilon\,\cG(u_t),
\]
which implies first $\cG(P_t^\varepsilon h)\le e^{C\varepsilon t}\cG(h)$
and then, after integration,
\begin{equation}\label{eq:staticeps}
 \cG(P_t^\varepsilon h)+c_*\int_0^t\cD(P_s^\varepsilon h)\dd s
 \le\cG(h)+C\varepsilon\int_0^t\cG(P_s^\varepsilon h)\dd s
 \le\big(1+C\varepsilon t\,e^{C\varepsilon t}\big)\cG(h).
\end{equation}
All constants here are independent of $\varepsilon$.

We now let $\varepsilon\to0$ for a fixed $t$. The strong convergence
\eqref{eq:visclimit}, the uniform bound on $\cG(P_t^\varepsilon h)$ and the
uniform bound on the time integral of $\cD(P_s^\varepsilon h)$ give weak
compactness of $\nabla_{x,v}P_t^\varepsilon h$ in $L^2(\mu)$ and of the
first and second derivatives appearing in $\cD$ in the time-integrated
space $L^2((0,t);L^2(\mu))$. The weak limits are the corresponding
distributional derivatives of $P_th$ and of $s\mapsto P_sh$, and weak lower
semicontinuity of the nonnegative quadratic forms on the left-hand side of
\eqref{eq:staticeps} yields \eqref{eq:static} for $h\in\Cc$, the error term
on the right-hand side vanishing in the limit. Finally, if $h\in H^1(\mu)$
we approximate it in $H^1(\mu)$ by $h_n\in\Cc$; then $\cG(h_n)\to\cG(h)$,
$P_sh_n\to P_sh$ in $L^2(\mu)$ uniformly in $s$, and the same lower
semicontinuity argument extends \eqref{eq:static} to $h$.
This completes the proof of Lemma~\ref{lem:static}. \qed

\begin{remark}
Inequality \eqref{eq:static} also shows that $(P_t)_{t\ge0}$ is a strongly
continuous semigroup on $H^1(\mu)$. Indeed, strong continuity in
$L^2(\mu)$ and the uniform bound $\cG(P_th)\le\cG(h)$ imply that
$P_th\to h$ weakly in $H^1(\mu)$ as $t\to0$, while the same bound and lower
semicontinuity imply $\cG(P_th)\to\cG(h)$; since $\cG$ is the square of a
Hilbert norm equivalent to that of $H^1(\mu)$, weak convergence together
with convergence of the norms gives strong convergence.
\end{remark}

\section{Proof of Lemma~\ref{lem:elliptic}}\label{apx:elliptic}

This appendix contains the proof of Lemma~\ref{lem:elliptic}. The material
is standard; it supplies, in the present setting, an approximation argument
of the kind which is omitted in \cite[Appendix~A.21]{Villani} and in
\cite{Zhang}.
We assume that Assumption~\ref{ass:WH} holds, and we fix
$0<\varepsilon\le1$;
constants denoted by $C_\varepsilon$ may depend on $\varepsilon$. Recall
that $\A_\varepsilon=\A+\varepsilon\Sx$ with
$\Sx=\sigma\Delta_x-\nabla L\cdot\nabla_x$, and that on test functions
\begin{equation}\label{eq:ellipticdisscore}
 -\Ree\ip{\A_\varepsilon h}h
 =\gamma\sigma\nm{\nabla_vh}^2+
 \varepsilon\sigma\nm{\nabla_xh}^2 ,
\end{equation}
by \eqref{eq:ibp} and by the symmetry of $\Sx$. In particular
$\A_\varepsilon$ is dissipative on $\Cc$ and closable, and we denote its
closure by the same symbol.

We shall use the following two multiplier bounds, which control the
unbounded coefficients $v$ and $\nabla L$ by one derivative in the
corresponding variable:
\begin{align}
 \nm{vh}^2&\le2d\sigma\nm h^2+4\sigma^2\nm{\nabla_vh}^2,
       &&h\in H^1_v(\mu),\label{eq:vmult}\\
 \nm{(\nabla L)h}^2
 &\le3d\sigma\sqrt{C_L}\nm h^2+
       (d\sigma\sqrt{C_L}+4\sigma^2)\nm{\nabla_xh}^2,
       &&h\in H^1_x(\mu).\label{eq:Lmult}
\end{align}
Let first $h\in\Cc$. Since
$v\,e^{-|v|^2/(2\sigma)}=-\sigma\nabla_ve^{-|v|^2/(2\sigma)}$, an
integration by parts in $v$ gives
\[
 \nm{vh}^2=d\sigma\nm h^2+
                  2\sigma\Ree\int\overline h\,v\cdot\nabla_vh\dd\mu
          \le d\sigma\nm h^2+\tfrac12\nm{vh}^2+2\sigma^2\nm{\nabla_vh}^2,
\]
where Young's inequality was used in the last step, and \eqref{eq:vmult}
follows after absorbing $\tfrac12\nm{vh}^2$. In the position variable the
analogous integration by parts, which rests on the identity
$(\nabla L)\,e^{-L/\sigma}=-\sigma\nabla_xe^{-L/\sigma}$, reads
\[
 \frac1\sigma\nm{(\nabla L)h}^2
   =\int\Delta L\,|h|^2\dd\mu
    +2\Ree\int\overline h\,\nabla L\cdot\nabla_xh\dd\mu .
\]
Since $|\Delta L|\le d\,\|\nabla^2L\|_{\rm op}$, the Cauchy--Schwarz
inequality and Assumption~\ref{ass:WH} bound the first integral by
$d\sqrt{C_L}(\nm h^2+\nm{\nabla_xh}^2)^{1/2}\nm h$, and therefore
\begin{align*}
 \frac1\sigma\nm{(\nabla L)h}^2
 &\le d\sqrt{C_L}\big(\nm h^2+\nm{\nabla_xh}\nm h\big)
                       +2\nm{(\nabla L)h}\nm{\nabla_xh}\\
 &\le\frac{3d\sqrt{C_L}}{2}\nm h^2+
       \frac{d\sqrt{C_L}}{2}\nm{\nabla_xh}^2
       +\frac1{2\sigma}\nm{(\nabla L)h}^2+
                                      2\sigma\nm{\nabla_xh}^2.
\end{align*}
Rearranging gives \eqref{eq:Lmult} for test functions. Both estimates extend
to $H^1_v(\mu)$ and $H^1_x(\mu)$ respectively by a density argument.

\subsection{Maximal dissipativity and the domain}

The formal adjoint of $\A_\varepsilon$ in $L^2(\mu)$ is
$\Lham+\gamma\LFD+\varepsilon\Sx=\J\A_\varepsilon\J$, because $\Sx$ is
symmetric and commutes with $\J$. Fix $\lambda>0$ and let $h\in L^2(\mu)$
be orthogonal to the range of $\lambda-\A_\varepsilon$ on $\Cc$. Then $h$
is a distributional solution of
$(\lambda-\Lham-\gamma\LFD-\varepsilon\Sx)h=0$, an equation whose
second-order part $\gamma\sigma\Delta_v+\varepsilon\sigma\Delta_x$ is
elliptic with smooth coefficients, and $h$ is therefore smooth. Let
$\chi_R$ be a Euclidean cutoff equal to one on the ball of radius $R$ and
zero outside the ball of radius $2R$, with $|\nabla\chi_R|\le C/R$. Testing
the equation against $\chi_R^2\overline h$ and taking real parts gives
\begin{equation}\label{eq:ellipticdeficiency}
 \begin{aligned}
 &\lambda\nm{\chi_Rh}^2+\gamma\sigma\nm{\chi_R\nabla_vh}^2
   +\varepsilon\sigma\nm{\chi_R\nabla_xh}^2\\
 &\quad=-2\gamma\sigma\Ree\ip{\chi_R\nabla_vh}{h\nabla_v\chi_R}
   -2\varepsilon\sigma\Ree\ip{\chi_R\nabla_xh}{h\nabla_x\chi_R}
   -\int|h|^2\chi_R\,\Lham\chi_R\dd\mu .
 \end{aligned}
\end{equation}
The first two terms on the right are bounded by Young's inequality by one
quarter of the two diffusion terms on the left plus $CR^{-2}\nm h^2$. For
the transport error we use
$|\Lham\chi_R|\le CR^{-1}(|v|+|\nabla L|)$ and the multiplier bounds
\eqref{eq:vmult}--\eqref{eq:Lmult}, applied to the test function
$\chi_Rh$:
\[
 \Big|\int|h|^2\chi_R\,\Lham\chi_R\dd\mu\Big|
 \le\frac CR\nm h\big(\nm{v\chi_Rh}+\nm{(\nabla L)\chi_Rh}\big)
               \le\frac CR\nm h\,\|\chi_Rh\|_{H^1(\mu)} .
\]
Since $\|\chi_Rh\|_{H^1(\mu)}\le\nm{\chi_R\nabla_{x,v}h}+C\nm h$, Young's
inequality absorbs the gradients into the two positive diffusion terms,
and what remains on the right-hand side is
$C_\varepsilon(R^{-1}+R^{-2})\nm h^2$, which tends to zero as
$R\to\infty$. Hence $h=0$, and the range of $\lambda-\A_\varepsilon$ on
$\Cc$ is dense. By \cite[Proposition~II.3.14(iv)]{EngelNagel}, the range of
$\lambda$ minus the closure of $\A_\varepsilon$ is the closure of this
range, hence all of $L^2(\mu)$, so that the closure of $\A_\varepsilon$ is
maximal dissipative. The argument of the proof of
Proposition~\ref{prop:adjoint} then applies without change: the adjoint of
$\A_\varepsilon$ is $\J\A_\varepsilon\J$, and the domain of
$\A_\varepsilon$ is the maximal domain.

For $f$ in the domain of $\A_\varepsilon$ and for functions $f_n\in\Cc$
such that $f_n\to f$ and $\A_\varepsilon f_n\to\A_\varepsilon f$ in
$L^2(\mu)$, \eqref{eq:ellipticdisscore} applied to
differences shows that $\nabla_{x,v}f_n$ converges in $L^2(\mu)$. Hence
$\Dom(\A_\varepsilon)\subset H^1(\mu)$, the identity
\eqref{eq:ellipticdisscore} extends to the whole domain, and its polarized
form is \eqref{eq:ellipticdiss}.

\subsection{Analyticity}

For $f,g\in H^1(\mu)$ we define
\[
 \mathfrak a_\varepsilon(f,g)=
 \gamma\sigma\ip{\nabla_vf}{\nabla_vg}
 +\varepsilon\sigma\ip{\nabla_xf}{\nabla_xg}
 +\ip{\nabla_xf}{vg}-\ip{\nabla_vf}{(\nabla L)g},
\]
This expression is linear in $f$ and conjugate linear in $g$, and it
satisfies $\mathfrak a_\varepsilon(f,g)=\ip{-\A_\varepsilon f}{g}$
for $f\in\Cc$ and $g\in H^1(\mu)$. By the multiplier inequalities
\eqref{eq:vmult}--\eqref{eq:Lmult} it is bounded on $H^1(\mu)$. The sum of
its last two terms has zero real part when $g=f$: this is the
antisymmetry of $\Lham$ for test functions, and it extends to $H^1(\mu)$ by
density and boundedness. The real part of $\mathfrak a_\varepsilon(f,f)$
is therefore
$\gamma\sigma\nm{\nabla_vf}^2+\varepsilon\sigma\nm{\nabla_xf}^2$, which
controls both gradients; after a scalar shift the norm induced by the real
part is equivalent to the norm of $H^1(\mu)$, and the imaginary part is
bounded by the square of that norm. By \cite[Chapters~VI and~IX]{Kato},
these two properties imply that there is a closed operator $T$ on
$L^2(\mu)$ such that
$\mathfrak a_\varepsilon(f,g)=\ip{Tf}{g}$ for $f$ in the domain of $T$ and
$g\in H^1(\mu)$, and $-T$ generates an analytic semigroup. The operator
$-T$ is closed and dissipative and extends $\A_\varepsilon$ on $\Cc$, hence it extends
the closure of $\A_\varepsilon$, and by maximality the two operators
coincide. This proves that $(P_t^\varepsilon)_{t\ge0}$ is an analytic
contraction semigroup, and completes the proof of part~(i) of
Lemma~\ref{lem:elliptic}.

\subsection{The second-derivative estimate}

Let $h\in\Cc$ and write
$\A_\varepsilon=(\gamma\LFD+\varepsilon\Sx)+(-\Lham)$ as the sum of its
symmetric and antisymmetric parts. The following two identities are global
second-derivative calculations of Bochner type, with all operators written
out:
\begin{align}
 &\nm{(\gamma\LFD+
             \varepsilon\sigma\Delta_x-\varepsilon\nabla L\cdot\nabla_x)h}^2
       \notag\\
 &\quad=\gamma^2\sigma^2\nm{\nabla_v^2h}^2+
       \gamma^2\sigma\nm{\nabla_vh}^2+
       2\varepsilon\gamma\sigma^2\nm{\nabla_x\nabla_vh}^2\notag\\
 &\qquad+\varepsilon^2\sigma^2\nm{\nabla_x^2h}^2+
       \varepsilon^2\sigma\Ree
                    \ip{(\nabla^2L)\nabla_xh}{\nabla_xh},
                                    \label{eq:ellipticbochner}\\
 &2\Ree\ip{
  (\gamma\LFD+
             \varepsilon\sigma\Delta_x-\varepsilon\nabla L\cdot\nabla_x)h}{
                (-\Lham)h}\notag\\
 &\quad=2\gamma\sigma\Ree\ip{\nabla_vh}{\nabla_xh}
       -2\varepsilon\sigma\Ree
                    \ip{\nabla_xh}{(\nabla^2L)\nabla_vh}.
                                    \label{eq:ellipticcross}
\end{align}
They are obtained by integrating by parts with \eqref{eq:ibp} and by using
the commutators $[\nabla_v,\LFD]=-\nabla_v$,
$[\nabla_x,\Sx]=-(\nabla^2L)\nabla_x$, $[\nabla_v,\Lham]=\nabla_x$ and
$[\nabla_x,\Lham]=-(\nabla^2L)\nabla_v$. Since
$\nm{\A_\varepsilon h}^2$ is the sum of \eqref{eq:ellipticbochner}, of
\eqref{eq:ellipticcross} and of the nonnegative term $\nm{\Lham h}^2$, we
may discard the latter. We then apply the relative bounds
\eqref{eq:relative} to the two terms containing the Hessian, and absorb the
resulting second derivatives by Young's inequality into the positive
second-order terms of \eqref{eq:ellipticbochner}. This gives
\[
 \nm{\A_\varepsilon h}^2
       \ge c_\varepsilon\nm{\nabla_{x,v}^2h}^2
                         -C_\varepsilon\nm{\nabla_{x,v}h}^2 ,
\]
and by \eqref{eq:ellipticdisscore} the last term is bounded by
$C_\varepsilon\nm{\A_\varepsilon h}\nm h$. This proves
\eqref{eq:ellipticgraph} for test functions. If $h$ belongs to the domain
of $\A_\varepsilon$ and the functions $h_n\in\Cc$ are such that
$h_n\to h$ and $\A_\varepsilon h_n\to\A_\varepsilon h$ in $L^2(\mu)$,
the estimate applied to differences shows that $h_n$ is a Cauchy sequence
in $H^2(\mu)$, so that $h\in H^2(\mu)$ and \eqref{eq:ellipticgraph} holds
for $h$. This proves part~(ii).

\subsection{Convergence as the regularization is removed}

For a test function $h$ and $\lambda>0$,
\[
 (\lambda-\A_\varepsilon)^{-1}(\lambda-\A)h-h
 =\varepsilon(\lambda-\A_\varepsilon)^{-1}
                 \Sx h\longrightarrow0
 \qquad(\varepsilon\to0),
\]
because the resolvents of the maximal dissipative operators
$\A_\varepsilon$ have norm at most $1/\lambda$. The functions
$(\lambda-\A)h$ with $h\in\Cc$ are dense in $L^2(\mu)$, since $\A$ is the
closure of its restriction to $\Cc$ and $\lambda-\A$ maps $\Dom(\A)$ onto
$L^2(\mu)$, and the uniform bound on
the resolvents therefore gives the strong convergence of
$(\lambda-\A_\varepsilon)^{-1}$ to $(\lambda-\A)^{-1}$. The
Trotter--Kato approximation theorem for contraction semigroups
\cite[Theorem~III.4.8]{EngelNagel} then proves
\eqref{eq:visclimit}, which is part~(iii).

\subsection{Regularity of the trajectories and the differentiation rule}

We finally prove part~(iv). The only property of the potential which is
needed is the pair of relative bounds \eqref{eq:relative}, which show that
multiplication by the Hessian of $L$ maps the gradient of a function in
$H^2(\mu)$ into $L^2(\mu)$.

\emph{Commutator bound.} Let $w\in H^2(\mu)$. The commutators of the first
derivatives with $\A$ and with $\Sx$, computed in the sense of
distributions, are
\begin{align*}
 [\nabla_v,\A]w&=-\nabla_xw-\gamma\nabla_vw,&
 [\nabla_x,\A]w&=(\nabla^2L)\nabla_vw,\\
 [\nabla_v,\Sx]w&=0,&
 [\nabla_x,\Sx]w&=-(\nabla^2L)\nabla_xw .
\end{align*}
By \eqref{eq:relative} they belong to $L^2(\mu)$, and for every first-order
derivative $\partial$ in the variables $(x,v)$ we have
\[
 \nm{[\partial,\A_\varepsilon]w}\le C\|w\|_{H^2(\mu)} .
\]

\emph{Bootstrap.} Suppose that $w\in\Dom(\A_\varepsilon)$, so that
$w\in H^2(\mu)$ by part~(ii), and that $\A_\varepsilon w\in H^1(\mu)$. For
a first-order derivative $\partial$, the distributional identity
\[
 \A_\varepsilon\partial w
   =\partial(\A_\varepsilon w)-[\partial,\A_\varepsilon]w
\]
has its right-hand side in $L^2(\mu)$, so that the maximal-domain
characterization puts $\partial w$ in $\Dom(\A_\varepsilon)$. Applying
\eqref{eq:ellipticgraph} to $\partial w$ and summing over the first-order
derivatives yields
\begin{equation}\label{eq:sobolevbootstrap}
 \|w\|_{H^3(\mu)}\le C_\varepsilon
            \big(\|\A_\varepsilon w\|_{H^1(\mu)}+\|w\|_{H^2(\mu)}\big)
 \le C_\varepsilon\big(\nm{\A_\varepsilon^2w}+\nm{\A_\varepsilon w}
            +\nm w\big),
\end{equation}
where the second inequality holds when in addition
$\A_\varepsilon w\in\Dom(\A_\varepsilon)$, and follows from
\eqref{eq:ellipticgraph} applied to $w$ and to $\A_\varepsilon w$.

\emph{Regularity of the trajectories.} Let $h\in\Cc$ and
$u_t=P_t^\varepsilon h$. All the generator powers
$\A_\varepsilon^ju_t=P_t^\varepsilon\A_\varepsilon^jh$, $j\ge0$, belong to
the domain of $\A_\varepsilon$ for every $t\ge0$, because each
$\A_\varepsilon^jh$ is again a test function; in particular they belong to
$H^2(\mu)$ by part~(ii). The bootstrap step therefore applies to $u_t$,
and more generally to every element $w$ of the linear span of the
functions $\A_\varepsilon^ju_t$ with $t\ge0$ and $j\ge0$, and it shows that
these functions belong to $H^3(\mu)$ and satisfy
\eqref{eq:sobolevbootstrap}. Moreover, for every $j\ge0$ the map
$t\mapsto\A_\varepsilon^ju_t$ is continuously differentiable on
$[0,\infty)$ with values in $L^2(\mu)$. Applying
\eqref{eq:sobolevbootstrap} to differences and to difference quotients of
the trajectory shows that $u\in C^1([0,\infty);H^3(\mu))$, with derivative
$\A_\varepsilon u_t$, including at $t=0$. This is the global regularity and
the time differentiability which are required before the energy
functionals of Appendix~\ref{apx:energy} can be differentiated.

\emph{The differentiation rule.} Let $\partial^\alpha$ and
$\partial^\beta$ be derivatives of order at most one. By the bootstrap
step, $\partial^\alpha u_t$ and $\partial^\beta u_t$ belong to
$\Dom(\A_\varepsilon)$, and by the previous step
$t\mapsto\partial^\alpha u_t$ is continuously differentiable with values
in $L^2(\mu)$, with derivative
\[
 \partial^\alpha\A_\varepsilon u_t
   =\A_\varepsilon\partial^\alpha u_t+[\partial^\alpha,\A_\varepsilon]u_t .
\]
Differentiating $\Ree\ip{\partial^\alpha u_t}{\partial^\beta u_t}$ and
applying \eqref{eq:ellipticdiss} to the pair
$(\partial^\alpha u_t,\partial^\beta u_t)$ gives \eqref{eq:star}. This
completes the proof of Lemma~\ref{lem:elliptic}. \qed

\section{Potentials of class $C^2$}\label{apx:C2}

This appendix contains the proof of the statement of
Remark~\ref{rem:C2}. Throughout it, $L\in C^2(\R^d)$ is bounded from below
and satisfies the pointwise condition \eqref{eq:VG} with a constant $C>0$,
and $\mu$ is the measure \eqref{eq:mu}. The differential expression $\A$
still maps $\Cc$ into $L^2(\mu)$, because its coefficients are continuous,
but Proposition~\ref{prop:maxdissipative} is no longer available, since it
requires a smooth potential. In its place we use a result of Conrad and
Grothaus \cite[Corollary~2.3 and Remark~2.4]{ConradGrothaus}: if the
potential is locally Lipschitz continuous and bounded from below, then the
operator $\A$, defined on $\Cc$, is closable in $L^2(\mu)$ and its closure
is maximal dissipative; the measure $\mu$ is not required to be finite. In
\cite{ConradGrothaus} this is stated for $\gamma=\sigma=1$ and for
real-valued functions. The general case reduces to it by the dilation
$x=(\sqrt\sigma/\gamma)\,x'$, $v=\sqrt\sigma\,v'$, which transforms
$\A/\gamma$ into the operator with $\gamma=\sigma=1$ and with the potential
$L(\sqrt\sigma\,x'/\gamma)/\sigma$, and the result carries over to
complex-valued functions because $\A$ has real coefficients. We denote the
closure again by $\A$. As in Section~\ref{sec:generator}, it generates a
strongly continuous semigroup of contractions $P_t$ on $L^2(\mu)$, and
$\Dom(\A^m)$ denotes the domain of the $m$th power of $\A$.

\begin{proposition}\label{prop:C2}
Let $L\in C^2(\R^d)$ be bounded from below and satisfy \eqref{eq:VG}.
Then all the conclusions of Theorem~\ref{thm:main} hold for the semigroup
$P_t$ generated by the closure of $\A$: for every $t>0$ and every integer
$m\ge1$ we have $P_tL^2(\mu)\subset\Dom(\A^m)$ and
$\A^mP_t=(\A P_{t/m})^m$, the operator $\A^mP_t$ maps $L^2(\mu)$ into
$H^1(\mu)$ for every $m\ge0$, and the estimates \eqref{eq:main-even} and
\eqref{eq:main-odd} hold, with constants $C_m$ and $D_m$ which depend only
on $m,d,\gamma,\sigma$ and on the constant $C$ of \eqref{eq:VG}.
\end{proposition}

The proof approximates $L$ by smooth potentials to which
Theorem~\ref{thm:main} applies with constants which do not depend on the
approximation, and it occupies the rest of this appendix.

\emph{Step 1: mollification.} Let $\rho\in C_c^\infty(\R^d)$ be a
nonnegative function which is supported in the unit ball and has integral
one, and let $\rho_\delta(y)=\delta^{-d}\rho(y/\delta)$ and
$L_\delta=L*\rho_\delta$ for $0<\delta\le1/(4C)$. The potential $L_\delta$
is smooth, and $L_\delta$ and $\nabla L_\delta$ converge to $L$ and
$\nabla L$ uniformly on compact sets as $\delta\to0$. We claim that
\begin{equation}\label{eq:VGdelta}
 \|\nabla^2L_\delta(x)\|_{\rm op}\le2C\big(1+|\nabla L_\delta(x)|\big),
 \qquad x\in\R^d .
\end{equation}
To see this, fix $x$, let $B$ be the closed ball of radius $\delta$ centred
at $x$, and let $M$ be the maximum of $|\nabla L|$ on $B$, which is
attained at some point $y_*\in B$. By \eqref{eq:VG} the Hessian of $L$ has
norm at most $C(1+M)$ on $B$. Since $B$ is convex and has diameter
$2\delta$, the mean value inequality gives, for every $y\in B$,
\[
 |\nabla L(y)-\nabla L(y_*)|\le2\delta\,C(1+M)\le\tfrac12(1+M) .
\]
The vector $\nabla L_\delta(x)$ is an average of the vectors $\nabla L(y)$,
$y\in B$, with the weight $\rho_\delta(x-y)$, so that the same bound holds
for $|\nabla L_\delta(x)-\nabla L(y_*)|$, and therefore
$|\nabla L_\delta(x)|\ge M-\tfrac12(1+M)$, that is,
$1+M\le2\big(1+|\nabla L_\delta(x)|\big)$. In the same way
$\nabla^2L_\delta(x)$ is an average of the matrices $\nabla^2L(y)$,
$y\in B$, whence
$\|\nabla^2L_\delta(x)\|_{\rm op}\le C(1+M)
\le2C\big(1+|\nabla L_\delta(x)|\big)$, which is \eqref{eq:VGdelta}.

By \eqref{eq:VGdelta} and \cite[Lemma~A.18]{Villani}, the smooth potential
$L_\delta$ satisfies Assumption~\ref{ass:WH} with a constant $C_L$ which
depends only on $d$, $\sigma$ and $C$, and in particular not on $\delta$.
We fix such a value of $C_L$. In the rest of the proof the constants of
Table~\ref{tab:constants} are evaluated at this value, so that they depend
only on $m,d,\gamma,\sigma$ and $C$.

\emph{Step 2: the approximating semigroups.} Let
$\dd\mu_\delta=e^{-(L_\delta+|v|^2/2)/\sigma}\dd x\dd v$ be the Gibbs
measure of $L_\delta$, let $\A_\delta$ be the generator associated with
$L_\delta$ in $L^2(\mu_\delta)$ by
Proposition~\ref{prop:maxdissipative}, and let $P^\delta_t$ be its
semigroup. By Step~1, Theorem~\ref{thm:main} and the estimates of
Section~\ref{sec:spectral} apply to $L_\delta$, and they give three bounds
which are uniform in $\delta$. For every $g\in L^2(\mu_\delta)$ and every
integer $m\ge0$,
\begin{align}
 \|\A_\delta P^\delta_t g\|_{\mu_\delta}
   &\le K\max\Big(1,\frac1t\Big)^2\|g\|_{\mu_\delta},&&t>0,
   \label{eq:unif1}\\
 \|\A_\delta^mP^\delta_tg\|_{\mu_\delta}
   &\le K_m\,t^{-m}\|g\|_{\mu_\delta},&&t\ge1,\label{eq:unif2}\\
 \|\nabla_{x,v}\A_\delta^mP^\delta_tg\|_{\mu_\delta}
   &\le M_m(t)\,\|g\|_{\mu_\delta},&&t>0,\label{eq:unif3}
\end{align}
where $M_m(t)^2=D_mt^{-(2m+1)}$ for $t\ge2$. The first bound is
\eqref{eq:APtau}, and the second is \eqref{eq:longnorm}, whose square is
\eqref{eq:main-even}. For $t\ge2$ the third bound is \eqref{eq:main-odd}.
For $0<t<2$ it follows from the argument of Section~\ref{sec:proofmain}:
with $s=\min\{1,t/2\}$ we have
$\A_\delta^mP^\delta_t=P^\delta_s\big(\A_\delta^mP^\delta_{t-s}\big)$,
the inner operator is bounded by \eqref{eq:powerinduction}, and
\eqref{eq:reverse-t} bounds the gradient of $P^\delta_sg$ by
$(2/(cs^3))^{1/2}\|g\|_{\mu_\delta}$; the resulting constant $M_m(t)$
depends only on $m$, $t$, $\gamma$, $\sigma$ and $C_L$.

The measures $\mu_\delta$ and $\mu$ differ, and we transport the
approximating objects to the fixed space $L^2(\mu)$. Since
$\dd\mu_\delta=e^{-(L_\delta-L)/\sigma}\dd\mu$, the multiplication operator
\[
 U_\delta g=e^{-(L_\delta-L)/(2\sigma)}\,g
\]
satisfies $\nm{U_\delta g}=\|g\|_{\mu_\delta}$, and it is a unitary map
from $L^2(\mu_\delta)$ onto $L^2(\mu)$. Consequently
$\widetilde P^\delta_t=U_\delta P^\delta_tU_\delta^{-1}$ is a strongly
continuous semigroup of contractions on $L^2(\mu)$, and its generator is
$\widetilde\A_\delta=U_\delta\A_\delta U_\delta^{-1}$, with the domain
$U_\delta\Dom(\A_\delta)$.

\emph{Step 3: convergence of the semigroups.} Let $h\in\Cc$. The function
$U_\delta^{-1}h=e^{(L_\delta-L)/(2\sigma)}h$ is twice continuously
differentiable and has compact support. It therefore belongs to the
maximal domain \eqref{eq:maxdomain} of $\A_\delta$, and $\A_\delta$ acts on
it as a differential operator. The factor $e^{(L_\delta-L)/(2\sigma)}$
depends on $x$ only, so that the only term of $\A_\delta$ which
differentiates it is $-v\cdot\nabla_x$, and we obtain
\[
 \widetilde\A_\delta h
 =\A_\delta h-\frac{v\cdot(\nabla L_\delta-\nabla L)}{2\sigma}\,h
 =\A h+(\nabla L_\delta-\nabla L)\cdot
       \Big(\nabla_vh-\frac{v}{2\sigma}\,h\Big),
\]
where the second equality uses
$\A_\delta h=\A h+(\nabla L_\delta-\nabla L)\cdot\nabla_vh$. Since
$\nabla L_\delta\to\nabla L$ uniformly on the support of $h$, this shows
that $\Cc\subset\Dom(\widetilde\A_\delta)$ and that
\begin{equation}\label{eq:genconv}
 \nm{\widetilde\A_\delta h-\A h}\longrightarrow0\qquad(\delta\to0),
 \qquad h\in\Cc .
\end{equation}
Since $\A$ is the closure of its restriction to $\Cc$, the Trotter--Kato
approximation theorem
\cite[Theorem~III.4.8]{EngelNagel} applies along every sequence
$\delta\to0$, and it gives
\begin{equation}\label{eq:sgconv}
 \nm{\widetilde P^\delta_th-P_th}\longrightarrow0\qquad(\delta\to0),
 \qquad h\in L^2(\mu),\quad t\ge0 .
\end{equation}

\emph{Step 4: the powers of the generator.} Fix $t>0$ and let
$\kappa=K\max(1,1/t)^2$. By \eqref{eq:unif1} and the unitarity of
$U_\delta$, the operator
$\widetilde\A_\delta\widetilde P^\delta_t
=U_\delta\big(\A_\delta P^\delta_t\big)U_\delta^{-1}$ is bounded on
$L^2(\mu)$, with norm at most $\kappa$. For $h\in\Cc$ the generator
$\widetilde\A_\delta$ commutes with its semigroup, and therefore
\[
 \nm{\widetilde\A_\delta\widetilde P^\delta_t h-P_t\A h}
 \le\nm{\widetilde\A_\delta h-\A h}
   +\nm{(\widetilde P^\delta_t-P_t)\A h}\longrightarrow0
\]
by \eqref{eq:genconv} and \eqref{eq:sgconv}. Since $P_t\A h=\A P_t h$,
we obtain $\nm{\A P_t h}\le\kappa\nm h$ for all $h\in\Cc$. Let now
$h\in L^2(\mu)$ and let $h_n\in\Cc$ converge to $h$. Then $P_t h_n$
converges to $P_t h$, and $\A P_t h_n$ is a Cauchy sequence by the
bound just proved. Since $\A$ is closed, this shows that
\begin{equation}\label{eq:C2smooth}
 P_t L^2(\mu)\subset\Dom(\A),\qquad
 \nm{\A P_t h}\le K\max\Big(1,\frac1t\Big)^2\nm h,
 \qquad h\in L^2(\mu),\quad t>0,
\end{equation}
and, since the operators $\widetilde\A_\delta\widetilde P^\delta_t$ are
bounded uniformly in $\delta$ and converge on the dense subspace $\Cc$,
that $\widetilde\A_\delta\widetilde P^\delta_t h\to\A P_t h$ for
every $h\in L^2(\mu)$.

The second part of the proof of Proposition~\ref{prop:smooth} uses only
that $\A P_s$ is a bounded operator for every $s>0$ and that $\A$ commutes
with the semigroup on $\Dom(\A)$. By \eqref{eq:C2smooth} it applies without
change, and it gives $P_tL^2(\mu)\subset\Dom(\A^m)$ and
$\A^mP_t=(\A P_{t/m})^m$ for $t>0$ and $m\ge1$. The same identity holds for
$\A_\delta$, and a product of operators which are bounded uniformly in
$\delta$ and converge strongly converges strongly to the product of the
limits. Hence, for every $h\in L^2(\mu)$, $t>0$ and $m\ge1$,
\begin{equation}\label{eq:powerconv}
 U_\delta\,\A_\delta^mP^\delta_t\,U_\delta^{-1}h
 =\big(\widetilde\A_\delta\widetilde P^\delta_{t/m}\big)^mh
 \longrightarrow(\A P_{t/m})^mh=\A^mP_th\qquad(\delta\to0),
\end{equation}
and for $m=0$ this convergence is \eqref{eq:sgconv}. Taking norms and using
\eqref{eq:unif2} together with
$\|U_\delta^{-1}h\|_{\mu_\delta}=\nm h$, we obtain
$\nm{\A^mP_th}\le K_mt^{-m}\nm h$ for $t\ge1$, which is
\eqref{eq:main-even} with $C_m=K_m^2$.

\emph{Step 5: the gradient.} Let $h\in L^2(\mu)$, $t>0$ and $m\ge0$, and
set $w_\delta=\A_\delta^mP^\delta_tU_\delta^{-1}h$ and $w=\A^mP_th$. By
\eqref{eq:unif3} the function $w_\delta$ belongs to $H^1(\mu_\delta)$ and
$\|\nabla_{x,v}w_\delta\|_{\mu_\delta}\le M_m(t)\nm h$, and by
\eqref{eq:powerconv} the functions $U_\delta w_\delta$ converge to $w$ in
$L^2(\mu)$. Let $\Phi=(\Phi^x,\Phi^v)$ be a continuously differentiable
vector field on $\R^{2d}$ with compact support and with values in
$\mathbb C^{2d}$, and let
$\nabla^*\Phi=\nabla_x^*\Phi^x+\nabla_v^*\Phi^v$, with the formal adjoints
\eqref{eq:adjoints}, so that $\nabla^*$ is the formal adjoint of
$\nabla_{x,v}$ in $L^2(\mu)$. The formal adjoint of $\nabla_{x,v}$ in
$L^2(\mu_\delta)$ maps $\Phi$ to
$\nabla^*\Phi+(\nabla L_\delta-\nabla L)\cdot\Phi^x/\sigma$. We test the
gradient of $w_\delta$ against the vector field
$U_\delta^{-1}\Phi=e^{(L_\delta-L)/(2\sigma)}\Phi$ in $L^2(\mu_\delta)$
and integrate by parts. The derivative of the exponential factor
contributes $-(\nabla L_\delta-\nabla L)\cdot\Phi^x/(2\sigma)$, and
returning to $L^2(\mu)$ we find
\[
 \ip{U_\delta\nabla_{x,v}w_\delta}{\Phi}
 =\Big\langle U_\delta w_\delta,\
   \nabla^*\Phi+\frac{(\nabla L_\delta-\nabla L)\cdot\Phi^x}{2\sigma}
  \Big\rangle_\mu .
\]
As $\delta\to0$ the right-hand side converges to $\ip w{\nabla^*\Phi}$,
because $U_\delta w_\delta\to w$ in $L^2(\mu)$ and
$\nabla L_\delta\to\nabla L$ uniformly on the support of $\Phi$. The
left-hand side is bounded in absolute value by
$\|\nabla_{x,v}w_\delta\|_{\mu_\delta}\nm\Phi\le M_m(t)\nm h\nm\Phi$.
Therefore
\[
 |\ip w{\nabla^*\Phi}|\le M_m(t)\,\nm h\,\nm\Phi
\]
for every such vector field $\Phi$. This means that the weak gradient of
$w$ belongs to $L^2(\mu)$ and that
$\nm{\nabla_{x,v}w}\le M_m(t)\nm h$. Thus $\A^mP_t$ maps $L^2(\mu)$ into
$H^1(\mu)$ for every $t>0$, and for $t\ge2$, where
$M_m(t)^2=D_mt^{-(2m+1)}$, this is the estimate \eqref{eq:main-odd}. The
proof of Proposition~\ref{prop:C2} is complete. \qed


\begin{thebibliography}{99}
\bibitem{FSW}
M. Fornasier, L. Sun, and R. Ward.
\emph{Large-time analysis of the Langevin dynamics for energies fulfilling
Polyak--{\L}ojasiewicz conditions}.
arXiv:2510.24925 (2025).
\url{https://arxiv.org/abs/2510.24925}.

\bibitem{CaoLuWang}
Y. Cao, J. Lu, and L. Wang.
\emph{On explicit $L^2$-convergence rate estimate for underdamped
Langevin dynamics}.
Archive for Rational Mechanics and Analysis \textbf{247}(5) (2023),
Article~90. \url{https://doi.org/10.1007/s00205-023-01922-4};
\url{https://arxiv.org/abs/1908.04746v7}.

\bibitem{HelfferNier}
B. Helffer and F. Nier.
\emph{Hypoelliptic Estimates and Spectral Theory for Fokker--Planck
Operators and Witten Laplacians}.
Lecture Notes in Mathematics, vol.~1862, Springer, 2005.
\url{https://doi.org/10.1007/b104762}.

\bibitem{NierLectures}
F. Nier.
\emph{Spectral analysis of Fokker--Planck and Kramers--Fokker--Planck
operators. Lecture summary}.
Lecture notes, January~9, 2012.
\url{https://www.math.univ-paris13.fr/~nier/recherche/exposes/resumeINRIA.pdf}.

\bibitem{HHS}
F. H\'erau, M. Hitrik, and J. Sj\"ostrand.
\emph{Tunnel effect and symmetries for Kramers--Fokker--Planck type operators}.
Journal of the Institute of Mathematics of Jussieu \textbf{10}(3) (2011),
567--634. \url{https://arxiv.org/abs/1007.0838}.

\bibitem{Hitrik}
M. Hitrik.
\emph{Tunnel effect and symmetries for non-selfadjoint operators}.
Journ\'ees \'equations aux d\'eriv\'ees partielles (2013), Expos\'e V, 1--12.
\url{https://doi.org/10.5802/jedp.101}.

\bibitem{NierChapter}
F. Nier.
\emph{Accurate estimates for the exponential decay of semigroups with
non-self-adjoint generators}.
Chapter~15 in O.~N. Kirillov and D.~E. Pelinovsky (eds.),
\emph{Nonlinear Physical Systems: Spectral Analysis, Stability and Bifurcations},
ISTE, 2014.
\url{https://onkirillov.com/wp-content/uploads/2015/03/ch15.pdf}.

\bibitem{ConradGrothaus}
F. Conrad and M. Grothaus.
\emph{Construction, ergodicity and rate of convergence of $N$-particle
Langevin dynamics with singular potentials}.
Journal of Evolution Equations \textbf{10} (2010), 623--662.
\url{https://doi.org/10.1007/s00028-010-0064-0}.

\bibitem{GrothausNonnenmacher}
M. Grothaus and A. Nonnenmacher.
\emph{Overdamped limit of generalized stochastic Hamiltonian systems for
singular interaction potentials}.
Journal of Evolution Equations \textbf{20} (2020), 577--605.
\url{https://doi.org/10.1007/s00028-019-00530-8};
\url{https://arxiv.org/abs/1809.06788}.
Statement numbers cited here refer to arXiv v1 (2018).

\bibitem{BreitenKunisch}
T. Breiten and K. Kunisch.
\emph{Improving the convergence rates for the kinetic Fokker--Planck
equation by optimal control}.
Preprint, \url{https://arxiv.org/abs/2205.01369}.
Statement numbers cited here refer to arXiv v2 (2025).

\bibitem{ABL}
P. Alphonse, J.-M. Bouclet, and M. L\'eautaud.
\emph{Short-time parametrix for the Fokker--Planck semigroup and
applications}.
Preprint, \url{https://arxiv.org/abs/2607.13814} (2026).

\bibitem{NierBoundary}
F. Nier.
\emph{Boundary Conditions and Subelliptic Estimates for Geometric
Kramers--Fokker--Planck Operators on Manifolds with Boundaries}.
Memoirs of the American Mathematical Society \textbf{252}, no.~1200 (2018).
\url{https://arxiv.org/abs/1309.5070}.
Section and proposition numbers cited here refer to arXiv v2 (2014).

\bibitem{HerauNier}
F. H\'erau and F. Nier.
\emph{Isotropic hypoellipticity and trend to equilibrium for the
Fokker--Planck equation with a high-degree potential}.
Archive for Rational Mechanics and Analysis \textbf{171} (2004), 151--218.
\url{https://doi.org/10.1007/s00205-003-0276-3}.

\bibitem{Herau}
F. H\'erau.
\emph{Short and long time behavior of the Fokker--Planck equation in a
confining potential and applications}.
Journal of Functional Analysis \textbf{244} (2007), 95--118.
\url{https://arxiv.org/abs/math/0501363}.

\bibitem{Villani}
C. Villani.
\emph{Hypocoercivity}.
Memoirs of the American Mathematical Society \textbf{202}, no.~950 (2009).
\url{https://arxiv.org/abs/math/0609050}.
Section and theorem numbers are those of the preprint.

\bibitem{BDMMS}
E. Bouin, J. Dolbeault, S. Mischler, C. Mouhot, and C. Schmeiser.
\emph{Hypocoercivity without confinement}.
Pure and Applied Analysis \textbf{2} (2020).
\url{https://doi.org/10.2140/paa.2020.2.203};
\url{https://arxiv.org/abs/1708.06180}.

\bibitem{BDS}
E. Bouin, J. Dolbeault, and C. Schmeiser.
\emph{Diffusion and kinetic transport with very weak confinement}.
Kinetic and Related Models \textbf{13}(2) (2020), 345--371.
\url{https://doi.org/10.3934/krm.2020012}.

\bibitem{Wang}
X.~P. Wang.
\emph{Large-time asymptotics of solutions to the Kramers--Fokker--Planck
equation with a short-range potential}.
Communications in Mathematical Physics \textbf{336} (2015), 1435--1471.
\url{https://doi.org/10.1007/s00220-014-2273-9};
\url{https://arxiv.org/abs/1403.0841}.

\bibitem{PWZ}
X. Pan, X.~P. Wang, and L. Zhu.
\emph{The Kramers--Fokker--Planck equation with a decaying potential in
$\R^n$, $n\ge4$}.
arXiv:2505.09270 (2025).
\url{https://arxiv.org/abs/2505.09270}.

\bibitem{Cao}
C. Cao.
\emph{The kinetic Fokker--Planck equation with weak confinement force}.
Communications in Mathematical Sciences \textbf{17}(8) (2019), 2281--2308.
\url{https://arxiv.org/abs/1801.10354}.

\bibitem{BSWW}
G. Brigati, G. Stoltz, A.~Q. Wang, and L. Wang.
\emph{Explicit convergence rates of underdamped Langevin dynamics under
weighted and weak Poincar\'e--Lions inequalities}.
arXiv:2407.16033v2, June 17, 2025.
\url{https://arxiv.org/abs/2407.16033v2}.

\bibitem{Zhang}
C. Zhang.
\emph{Hypocoercivity and global hypoellipticity for the kinetic
Fokker--Planck equation in $H^k$ spaces}.
Kinetic and Related Models \textbf{17}(3) (2024), 394--435.
\url{https://doi.org/10.3934/krm.2023027};
\url{https://arxiv.org/abs/2012.06253}.
All numbered references to this work use arXiv:2012.06253v1.

\bibitem{BattyDuyckaerts}
C.~J.~K. Batty and T. Duyckaerts.
\emph{Non-uniform stability for bounded semi-groups on Banach spaces}.
Journal of Evolution Equations \textbf{8} (2008), 765--780.

\bibitem{BorichevTomilov}
A. Borichev and Y. Tomilov.
\emph{Optimal polynomial decay of functions and operator semigroups}.
Mathematische Annalen \textbf{347} (2010), 455--478.
\url{https://doi.org/10.1007/s00208-009-0439-0}.

\bibitem{RSS}
J. Rozendaal, D. Seifert, and R. Stahn.
\emph{Optimal rates of decay for operator semigroups on Hilbert spaces}.
Advances in Mathematics \textbf{346} (2019), 359--388.
\url{https://doi.org/10.1016/j.aim.2019.02.007}.

\bibitem{CST}
R. Chill, D. Seifert, and Y. Tomilov.
\emph{Semi-uniform stability of operator semigroups and energy decay of
damped waves}.
Philosophical Transactions of the Royal Society A \textbf{378} (2020).
\url{https://arxiv.org/abs/2007.04711}.

\bibitem{BCT}
C.~J.~K. Batty, R. Chill, and Y. Tomilov.
\emph{Fine scales of decay of operator semigroups}.
Journal of the European Mathematical Society \textbf{18}(4) (2016), 853--929.
\url{https://doi.org/10.4171/JEMS/605};
\url{https://arxiv.org/abs/1305.5365}.
Theorem numbers refer to arXiv v2.

\bibitem{BobkovLedoux}
S.~G. Bobkov and M. Ledoux.
\emph{Weighted Poincar\'e-type inequalities for Cauchy and other convex measures}.
The Annals of Probability \textbf{37}(2) (2009), 403--427.
\url{https://doi.org/10.1214/08-AOP407};
\url{https://arxiv.org/abs/0906.1651}.

\bibitem{Ahlfors}
L.~V. Ahlfors.
\emph{Complex Analysis}. Third edition, McGraw--Hill, 1979.

\bibitem{Phillips}
R.~S. Phillips.
\emph{Dissipative operators and hyperbolic systems of partial differential
equations}.
Transactions of the American Mathematical Society \textbf{90} (1959),
193--254.

\bibitem{EngelNagel}
K.-J. Engel and R. Nagel.
\emph{One-Parameter Semigroups for Linear Evolution Equations}.
Graduate Texts in Mathematics, vol.~194, Springer, New York, 2000.
\url{https://doi.org/10.1007/b97696}.

\bibitem{Kato}
T. Kato.
\emph{Perturbation Theory for Linear Operators}.
Second edition, Grundlehren der mathematischen Wissenschaften, vol.~132,
Springer, Berlin, 1976.
\end{thebibliography}
\end{document}